\documentclass{amsart}

\def\makeatother{\catcode64=\active}
\def\mymessage#1{\immediate\write16{macros: #1}}

\usepackage{pstricks}
\usepackage{pst-node}
\usepackage{amscd}

\usepackage{amssymb}

\mymessage{Using times for text}
\usepackage{times}

\let\boldit\boldsymbol

\newcommand{\UseRsfsAsCalli}{
\DeclareFontFamily{U}{rsfs}{}
\DeclareFontShape{U}{rsfs}{m}{n}{%
   <5>rsfs5%
   <6>rsfs10%
   <7>rsfs7%
   <8>rsfs10%
   <9>rsfs10%
   <10>rsfs10%
   <11>rsfs10%
   <12>rsfs10%
   <14>rsfs10%
   <17>rsfs10%
   <20>rsfs10%
   <25>rsfs10}{}
\DeclareMathAlphabet{\callig}{U}{rsfs}{m}{n}
\newcommand{\calli}[1]{{\callig ##1\/}}
\mymessage{Using rsfs as calli}
}

\UseRsfsAsCalli

\DeclareFontFamily{U}{wncy}{}
\DeclareFontShape{U}{wncy}{m}{n}{%
   <5>wncyr5%
   <6>wncyr6%
   <7>wncyr7%
   <8>wncyr8%
   <9>wncyr9%
   <10>wncyr10%
   <11>wncyr10%
   <12>wncyr6%
   <14>wncyr7%
   <17>wncyr8%
   <20>wncyr10%
   <25>wncyr10}{}
\DeclareMathAlphabet{\cyrille}{U}{wncy}{m}{n}
\mymessage{Defining cyrille}

\newbox\dummybox
\def\mysubscripts{
\mymessage{WARNING: moving subscripts lower}
\setbox\dummybox\hbox{$$\fontdimen16\textfont2=2.8pt$$
$$\fontdimen16\scriptfont2=1.9pt$$}}
\mysubscripts

\newcommand{\clap}[1]{\hbox to 0pt{\hss #1\hss}}
\newcommand{\mathclap}[1]{{\mathchoice{\clap{$\displaystyle #1$}}
{\clap{$#1$}}
{\clap{$\scriptstyle #1$}}
{\clap{$\scriptscriptstyle #1$}}}}

\newcommand{\downiso}{\clap{$\left\downarrow\vphantom{\Bigl(}\right.$}
 \rlap{\hbox{}\raise0.3ex\hbox{$\wr$}}}

\newcommand{\upsurj}{\clap{$\left\uparrow\vphantom{\Bigl(}\right.$}
 \raise 0.1ex\clap{$\left\uparrow\vphantom{\bigl(}\right.$}}
\newcommand{\longsearrow}{\lower 1.4ex\hbox{\begin{picture}(18,18)(0,0)
\put(0,18){\vector(1,-1){18}}
\end{picture}}}
\newcommand{\longseearrow}{\lower 1.3ex\hbox{\begin{picture}(30,15)(0,0)
\put(0,15){\vector(2,-1){30}}
\end{picture}}}

\newbox\isobox
\newdimen\isodim
\newdimen\isoscriptdim
\newdimen\isoscriptscriptdim
\newcommand{\iso}{\mathrel{
\setbox\isobox\hbox{$\longrightarrow$}
\isodim=\wd\isobox
\setbox\isobox\hbox{$\scriptstyle\longrightarrow$}
\isoscriptdim=\wd\isobox
\setbox\isobox\hbox{$\scriptscriptstyle\longrightarrow$}
\isoscriptscriptdim=\wd\isobox
\mathchoice
{\hbox to\isodim{\hfil\lower 0.2ex\clap{$\widetilde{}$}%
\clap{$\longrightarrow$}\hfil}}%
{\hbox to\isodim{\hfil\lower 0.2ex\clap{$\widetilde{}$}%
\clap{$\longrightarrow$}\hfil}}%
{\hbox to\isoscriptdim{\hfil\lower 0.50ex%
\clap{$\scriptstyle\widetilde{}$}%
\clap{$\scriptstyle\longrightarrow$}\hfil}}%
{\hbox to\isoscriptscriptdim{\hfil\lower 0.60ex%
\clap{$\tilde{}$}%
\clap{$\scriptscriptstyle\longrightarrow$}\hfil}}}}

\catcode`@=11
\atdef@ surj{\ampersand@
  \ifCD@ \global\bigaw@\minCDarrowwidth \else \global\bigaw@\minaw@ \fi
 \ifCD@\enskip\fi
   \mathrel{\mathop{\hbox to\bigaw@{\rightarrowfill@\displaystyle
     \kern-2ex{$\to$}}}}%
 \ifCD@\enskip\fi \ampersand@}
\atdef@ i#1i#2i{\ampersand@
  \ifCD@ \global\bigaw@\minCDarrowwidth \else 
  \global\bigaw@\minaw@\fi
  \setbox\isobox\hbox{$\widetilde{\hbox{}}$}%
  \ifdim\wd\isobox>\bigaw@\global\bigaw@\wd\isobox\fi
  \setboxz@h{$\m@th\scriptstyle\;{#1}\;\;$}%
  \ifdim\wdz@>\bigaw@\global\bigaw@\wdz@\fi
  \@ifnotempty{#2}{\setbox\@ne\hbox{$\m@th\scriptstyle\;{#2}\;\;$}%
    \ifdim\wd\@ne>\bigaw@\global\bigaw@\wd\@ne\fi}%
  \ifCD@\enskip\fi
  \isodim\bigaw@
  \divide\isodim by 2
  \mathrel{\mathop{\hbox to\bigaw@{\hfil
        \clap{\hbox to \bigaw@{\rightarrowfill@\displaystyle}}%
        \lower 0.5ex\clap{$\widetilde{\hbox to\isodim{\hfil}}$}%
        \hfil}}%
    \limits^{#1}\@ifnotempty{#2}{_{#2}}}%
  \ifCD@\enskip\fi \ampersand@}
\makeatother

\CDat

\catcode`@=11
\atdef@ @{@}
\makeatother

\newbox\fordim
\newdimen\notdividim
\newdimen\notdiviscriptdim
\newdimen\notdiviscriptscriptdim
\setbox\fordim\hbox{$/$}
\notdividim=\wd\fordim
\setbox\fordim\hbox{$\scriptstyle/$}
\notdiviscriptdim=\wd\fordim
\setbox\fordim\hbox{$\scriptscriptstyle/$}
\notdiviscriptscriptdim=\wd\fordim
\newcommand{\notdivise}{%
\mathrel{\mathchoice%
{{\hbox to\notdividim{\hfil\clap{$/$}\clap{$|$}\hfil}}}%
{{\hbox to\notdividim{\hfil\clap{$/$}\clap{$|$}\hfil}}}%
{{\hbox to\notdiviscriptdim{\hfil\clap{$\scriptstyle/$}%
\clap{$\scriptstyle|$}\hfil}}}%
{{\hbox to\notdiviscriptscriptdim{\hfil\clap{$\scriptscriptstyle/$}%
\clap{$\scriptscriptstyle|$}\hfil}}}}}
\newcommand{\divise}{\mathrel
{\mathchoice{\hbox to -0.2em{\hss$\displaystyle|$\hss}}%
{\hbox to -0.2em{\hss$\textstyle|$\hss}}%
{\hbox to 0.05em{\hss$\scriptstyle|$\hss}}%
{\hbox to 0.05em{\hss$\scriptscriptstyle|$\hss}}}}

\catcode`@=11
\numberwithin{equation}{subsection}
\@addtoreset{equation}{section}
\makeatother

\catcode`@=11
\renewcommand\thesubsection{\ifnum\c@subsection=0\relax\thesection\else
\thesection.\arabic{subsection}\fi}
\makeatother

\newenvironment{scriptarray}{%
\setbox\strutbox\hbox{\vphantom{$\scriptstyle ($}}

\def\crr{\\\scriptstyle}
\begin{array}{c}\scriptstyle}{\end{array}}

\newdimen\graphsize
\DeclareMathOperator{\Det}{Det}
\DeclareMathOperator{\Trace}{Tr}
\DeclareMathOperator{\Pic}{Pic}

\DeclareMathOperator{\Spec}{Spec}

\DeclareMathOperator{\Fr}{Fr}
\DeclareMathOperator{\re}{Re}
\DeclareMathOperator{\rk}{rk}
\DeclareMathOperator{\Br}{Br}
\DeclareMathOperator{\Gal}{Gal}

\DeclareMathOperator{\Val}{Val}
\DeclareMathOperator{\Vol}{Vol}
\DeclareMathOperator{\Area}{Area}

\newcommand{\varleq}{\leqslant}
\newcommand{\vargeq}{\geqslant}
\renewcommand{\leq}{\leqslant}

\mymessage{Redefining leq and geq}

\newtheorem{theo}{theorem}[section]
\newtheorem{lem}[theo]{Lemma}
\newtheorem{prop}[theo]{Proposition}

\theoremstyle{definition}
\newtheorem{defi}[theo]{Definition}

\theoremstyle{remark}
\newtheorem{rem}[theo]{Remark}

\newtheorem{notas}[theo]{Notations}
\catcode`@=11
\@addtoreset{theo}{section}
\makeatother

\newcommand{\noqed}{\renewcommand{\qed}{}}

\newcommand{\cardinal}{\sharp}
\newcommand{\dual}{^\vee}

\newcommand{\ZZ}{{\mathbf Z}}
\newcommand{\QQ}{{\mathbf Q}}
\newcommand{\RR}{{\mathbf R}}
\newcommand{\CC}{{\mathbf C}}
\newcommand{\Ov}{{\calli O}_v}
\newcommand{\Adeles}{{\boldsymbol A}}
\newcommand{\FF}{{\mathbf F}}
\newcommand{\Ffp}{\FF_{\mathfrak p}}
\newcommand{\Fp}{\FF_{p}}

\newcommand{\Ceff}{\Lambda_{\text{\rm eff}}}
\newcommand{\Vbar}{{\overline V}}
\newcommand{\antican}{{\omega_V^{-1}}}

\newcommand{\normof}[1]{\Vert #1\Vert_v}
\newcommand{\metric}{\normof\cdot}
\newcommand{\height}{H}
\newcommand{\nUH}{N_{U,H}}

\newcommand{\Haar}[1]{{\text{d}#1\,}}
\newcommand{\omegaH}{\boldsymbol \omega_{\height}}
\newcommand{\omegaHof}[1]{\boldsymbol \omega_{\height,#1}}
\newcommand{\omegaHv}{\boldsymbol \omega_{\height,v}}

\newcommand{\omegaHp}{\boldsymbol \omega_{\height,p}}
\newcommand{\tauH}{\tau_H}
\newcommand{\thetaH}{\boldsymbol \theta_{\height}}

\catcode`@=11
\newtoks\temporary
\long\def\gaddtomacro#1#2{%
\begingroup
\temporary\expandafter{#1#2}%
\xdef#1{\the\temporary}
\endgroup}
\newcount\linenumber
\newcount\tobedone
\newcount\tablesize
\def\tablename{default}
\def\thelinenumber{\ifcase\linenumber zero\or one\or two\or
three\or four\or five\or six\or seven\or eight\or nine\or ten\or
eleven\or twelve\or thirteen\or fourteen\or fifteen\or sixteen\or
seventeen\or eighteen\or nineteen\or twenty\or twentyone\or
twentytwo\or twentythree\or twentyfour\or twentyfive\or
twentysix\or twentyseven\or twentyeight\or twentynine\fi}
\def\clearone{%
\expandafter\def\csname\tablename\thelinenumber\endcsname{}%
\advance\linenumber by 1}
\def\newone#1{%
\expandafter\gaddtomacro\csname\tablename
\thelinenumber\endcsname{#1}%
\advance\linenumber by 1}
\def\expandone{\csname\tablename\thelinenumber\endcsname
\global\advance\linenumber by 1}
\def\repeatgobble#1{\advance\tobedone by -1
\ifnum\the\tobedone>0\def\mynext{\repeatgobble}\else
\def\mynext{\relax}\fi\mynext}
\def\repeatclear{\clearone\advance\tobedone by -1
\ifnum\the\tobedone>0\def\mynext{\repeatclear}\else
\def\mynext{\relax}\fi\mynext}
\def\repeatnewone#1{\newone{#1}%
\advance\tobedone by -1
\ifnum\the\tobedone>0\def\mynext{\repeatnewone}\else
\def\mynext{\relax}\fi\mynext}
\def\repeatconstant#1{\newone{#1}%
\advance\tobedone by -1
\ifnum\the\tobedone>0\def\mynext{\repeatconstant{#1}}\else
\def\mynext{\relax}\fi\mynext}
\def\repeatexpandone{\expandone
\global\advance\tobedone by -1
\ifnum\the\tobedone>0%
\def\mynext{\crcr\tablerule\repeatexpandone}\else
\def\mynext{\crcr\tablerule}\fi\mynext}
\def\dogobble{\tobedone\tablesize
\linenumber 1
\repeatgobble}
\def\doclear{\tobedone\tablesize
\linenumber 1
\repeatclear}
\def\doconstant#1{\tobedone\tablesize
\linenumber 1
\repeatconstant{#1}}
\def\doaddcolumn{\tobedone\tablesize
\linenumber 1
\repeatnewone}
\def\docolumns{\global\tobedone\tablesize
\global\linenumber 1
\repeatexpandone}
\def\doinit#1#2{\def\tablename{#1}%
\tablesize#2
\doclear}
\def\tablerule{\noalign{\hrule}}
\def\linebegin{\global\def\intercolumn{\global\def\intercolumn{&}}}
\newcommand{\bigstrut}{\vphantom{$\Bigl)$}}
\def\createnewtable#1#2{%
\expandafter\def\csname for#1\endcsname{%
\doinit{#1}{#2}%
\doconstant{\linebegin}
\expandafter\def\csname add#1column\endcsname{%
\def\tablename{#1}%
\tablesize #2
\expandafter\def\csname First#1column\endcsname{%
\def\tablename{#1}%
\tablesize #2
\dogobble}
\doconstant{\intercolumn}
\doaddcolumn}
\expandafter\def\csname First#1column\endcsname{%
\csname add#1column\endcsname}
\expandafter\def\csname show#1columns\endcsname{%
\def\tablename{#1}%
\tablesize #2
\[\hbox to\hsize{\hss\vbox{%
\halign{\vrule\hskip1em\bigstrut
$########$\hfil\hskip1em\vrule
\tabskip=0pt%
&&\hskip0.5em\bigstrut$########$\hfil
\hskip0.5em\vrule\crcr\tablerule
\docolumns
}}\hss}\]
}}}
\createnewtable{HB}{14}
\createnewtable{CTKS}{17}
\newbox\trash
\newif\ifShowInfo
\ShowInfofalse
{\ifShowInfo\else\egroup\endgroup\fi}

\makeatother

\author{Emmanuel Peyre}
\address{Institut Fourier\\
UFR de Math\'ematiques, UMR 5582\\
Universit\'e de Grenoble I et CNRS\\
BP 74\\ 38402 Saint-Martin d'H\`eres CEDEX\\ France}
\urladdr{http://www-fourier.ujf-grenoble.fr/\~{}peyre}
\email{Emmanuel.Peyre@@ujf-grenoble.fr}
 
\author{Yuri Tschinkel}
\address{Department of Mathematics\\
Princeton University\\
Washington Rd.\!\\
Princeton, NJ 08544-1000\\
U.S.A.}
\email{ytschink@@math.princeton.edu}

\title[Tamagawa numbers of cubic
surfaces]{Tamagawa numbers of diagonal
cubic surfaces\\ of higher rank}
\date{\today}

\begin{document}
\begin{abstract}
We consider diagonal cubic surfaces defined
by an equation of the form
\[ax^3+by^3+cz^3+dt^3 = 0.\]
Numerically, one can find all rational points of
height $\varleq B$ for $B$ in the range of up to  $10^5$,
thanks to a program due to D. J. Bernstein.
On the other hand, there are precise conjectures concerning
the constants in the asymptotics of rational points of bounded
height due to Manin, Batyrev and the authors.
Changing the coefficients one can obtain cubic surfaces
with rank of the Picard group
varying between 1 and 4. We check that numerical data
are compatible with the above conjectures.
In a previous paper we considered cubic surfaces
with Picard groups of rank one with or without Brauer-Manin
obstruction to weak approximation. In this paper, we test the conjectures
for diagonal cubic surfaces with Picard groups of higher rank.
\end{abstract}

\maketitle
\tableofcontents

\section{Introduction}

This paper is devoted to numerical tests of a refined version
of a conjecture of Manin about the number of points of bounded
height on Fano varieties (see \cite{batyrevmanin:hauteur},
\cite{fmt:fano}, \cite{peyre:fano}, or \cite{batyrevtschinkel:tamagawa}
for a description of the conjectures). The choice of diagonal cubic 
surfaces to test these conjectures was motivated by the work of Heath-Brown
\cite{heathbrown:density} in which he treated the cases
\[X^3+Y^3+Z^3+aT^3=0\]
for $a=2$ or $3$. The results he obtained were used as a benchmark
for the subsequent attempts to interpret the asymptotic constants
(see, in particular, \cite{swinnertondyer:cubic},
\cite{peyre:fano} and \cite{peyretschinkel:numericone}).
\par
More precisely, we consider a diagonal cubic surface 
$V\subset\mathbf P^3_{\QQ}$
given by an equation of the form
\[aX^3+bY^3+cZ^3+dT^3=0.\]
Let $H$ be the height function on 
$\mathbf P^3(\QQ)$ defined by the formula:
for any $Q=(x:y:z:t)$ in $\mathbf P^3(\QQ)$, one has
\[
H(Q)=\max\{\vert x\vert,\vert y\vert,\vert z\vert,\vert t\vert\}
\text{ if }\begin{cases}
(x,y,z,t)\in\ZZ^4, \\
\gcd(x,y,z,t)=1.
\end{cases}\]
Let $U$ be the complement in $V$ to the 27 lines. We are interested
in the asymptotic behavior of the cardinal
\[\nUH(B)=\cardinal\{\,Q\in U(\QQ)\mid H(Q)\varleq B\,\}\]
as $B$ goes to infinity.
\par
Assume that $V(\QQ)$ is Zariski dense, which by a result of Segre 
(see \cite[\S 29,\S 30]{manin:cubic}) is equivalent
to $V(\QQ)\neq\emptyset$. It is expected that
\[\nUH(B)=BP(\log(B))+o(B)\]
as $B$ goes to $+\infty$, where $P$ is a polynomial of degree
$\rk\Pic(V)-1$, with leading coefficient $\thetaH(V)$. 
This constant has a conjectural description. 
The goal is to compute $\thetaH(V)$ explicitly in the examples at hand
and to compare it with numerical data. Our previous paper 
\cite{peyretschinkel:numericone} was devoted to surfaces
with Picard groups of rank one with or without Brauer-Manin obstruction
to weak approximation. In this paper, we consider examples
with Picard groups of higher rank. Note 
that in these examples the error
term 
\[(\nUH(B)-\thetaH(V)B(\log B)^{\rk\Pic(V)-1})/B(\log B)^{\rk\Pic(V)-1})\]
is expected to decrease more slowly.
Indeed, if $\rk\Pic(V)=1$ this error term
is expected to decrease as $1/B^{\epsilon }$ for some $\epsilon>0$,
whereas for higher ranks it should be comparable to $1/\log B$.
However we observe a very good accordance, which is even more striking
if we take into account that a polynomial $P$ of degree
$\rk\Pic(V)-1$ should appear in the asymptotics and
use a rather straightforward statistical formula
to estimate the leading coefficient of this polynomial 
from numerical data.
\par
The paper is organized as follows: 
in section \ref{section:constant} we define $\thetaH(V)$. Section \ref{section:lines}
contains the description of the Galois action on the geometric Picard
group $\Pic(\Vbar)$. In section~\ref{section:euler} we compute
the Euler product corresponding to good reduction places. 
In section \ref{section:bad} we explain how to compute
the local densities at the places of bad reduction.
In section \ref{section:alpha} we determine in each case the value
of the geometric constant $\alpha(V)$ defined in \S\ref{section:constant}. 
Section \ref{section:stats} is devoted to the description
of statistical tools we used to analyze the numerical data. In section
\ref{section:results} we present the results.

\section{Description of the conjectural constant}
\label{section:constant}

In this section we give a short description of the conjectural
asymptotic constant for heights defined
by an adelic metrization of the anticanonical line bundle 
(see \cite{peyre:fano} for more details and \cite{batyrevtschinkel:tamagawa}
for a discussion in a more general setting).

\begin{notas}
For any field $E$, we denote by $\overline E$ an algebraic closure
of $E$. If $X$ is a variety over $E$, then $X(E)$ denotes the set 
of rational points of $X$ and 
$\overline X$ the product
\hbox{$X\times_{\Spec(E)}\Spec\overline E$}. The cohomological Brauer
group $\Br (X)$ is defined as the \'etale cohomology group
$H^2_{\text{\'et}}(X,\mathbf G_m)$. For any $A$ in $\Br (X)$,
any extension $E'$ of $E$ and any $P$ in $V(E')$, we denote
by $A(P)$ the evaluation of $A$ at $P$.
\par
For a number field $F$ we denote by $\Val(F)$ the set of places
of $F$ and by $\Val_f(F)$ the set of finite places. 
The absolute discriminant of $F$ is denoted by $d_F$.
For any place $v$
of $F$, let $F_v$ be the $v$-adic completion of $F$. If $v$ is finite,
then $\Ov$ is the ring of $v$-adic integers and $\mathbf F_v$
the residue field. By global class field theory we have an exact sequence
\begin{equation}
\label{equ:constant:classfield}
0\to\Br(F)\to\!\!\!\bigoplus_{v\in\Val(F)}\!\!\!\Br(F_v)@>\sum\text{inv}_v>>
\mathbf\QQ/\mathbf\ZZ\to0.
\end{equation}
In this paragraph, $V$ is a smooth projective
geometrically integral variety over a number field $F$
satisfying the conditions:
\begin{itemize}
\item[(i)]$H^i(V,\calli O_V)=0$ for $i=1$ or $2$,
\item[(ii)]$\Pic(\overline V)$ has no torsion,
\item[(iii)]the anticanonical line bundle $\antican$
belongs to the interior of the cone of classes of effective
divisors $\Ceff(V)\subset\Pic(V)\otimes_\ZZ\RR$.
\end{itemize}
\par
The adelic space $V(\Adeles_F)$ of $V$ coincides with
the product $\prod_{v\in\Val(F)}V(F_v)$.
By \cite[lemma 1]{colliotthelene:hasse}, for
any class $A$ in $\Br(V)$, one has a map $\rho_A$
defined as the composition
\[\begin{array}{rcccc}
V(\Adeles_F)&\to&\bigoplus_{v\in\Val(F)}\Br(F_v)&@>\sum\text{inv}_v>>
\QQ/\ZZ\\
(P_v)_{v\in\Val(F)}&\mapsto&(A(P_v))_{v\in\Val(F)}.
\end{array}\]
Then one defines
\[V(\Adeles_F)^{\Br}=\bigcap_{A\in\Br(V)}\ker(\rho_A)\subset V(\Adeles_F).\]
By the exact sequence \eqref{equ:constant:classfield}, one has the inclusion
\[\overline{V(F)}\subset V(\Adeles_F)^{\Br}\]
where $\overline{V(F)}$ denotes the topological closure 
of the set of rational points. Conjecturally both sets coincide
for cubic surfaces.
(See also the text of Swinnerton-Dyer in this volume).
The {\em Brauer-Manin obstruction to weak approximation}, described
by Manin in \cite{manin:brauer} and by Colliot-Th\'el\`ene and Sansuc
in \cite{cts:predescente2}, is the condition
\[V(\Adeles_F)^{\Br}\neq V(\Adeles_F).\]
\par
Let us assume that the height $\height$ on $V$ is defined
by an adelic metric $(\metric)_{v\in\Val(F)}$ on $\antican$.
By definition, this means that we consider
$\antican$ as a line bundle, that the functions $\metric$
are $v$-adically continuous metrics on $\antican(F_v)$
which for almost all places $v$ are given by a smooth model of $V$,
and that the height of a rational point $x$ of $V$ is given
by the formula
\[\forall y\in\antican(x)-\{0\},\quad H(x)=\prod_{v\in\Val(F)}
\normof y^{-1}\]
where $\antican(x)$ is the fiber of $\antican$ at $x$.
\par
If $v\in \Val(F)$ the Haar measure $\Haar {x_v}$
on $F_v$ is normalized as follows:
\begin{itemize}
\item[-]$\int_{\Ov}\Haar{x_v}=1$ if $v$ is finite,
\item[-]$\Haar{x_v}$ is the usual Lebesgue measure if $F_v\iso\RR$,
\item[-]$\Haar{x_v}=\Haar{z}\Haar{\overline z}=2\Haar x\Haar y$
if $F_v\iso\CC$.
\end{itemize}
The metric $\metric$ defines a measure $\omegaHv$ on the
locally compact space $V(F_v)$. In local $v$-adic analytic coordinates
$x_{1,v}\dots x_{n,v}$ on $V(F_v)$
this measure is given by the formula
\[\omegaHv=\left\Vert\frac{\partial}{\partial x_{1,v}}\wedge\dots\wedge
\frac{\partial}{\partial x_{n,v}}\right\Vert_v\Haar{x_{1,v}}\dots
\Haar{x_{n,v}}.\]
\par
If $M$ is
a discrete representation of  $\Gal(\overline F/F)$ over $\QQ$,
then for any finite place $\mathfrak p$ of $F$,
the local term of the corresponding Artin 
$L$-function is defined as follows:
we choose an algebraic closure $\overline F_{\mathfrak p}$
of $F_\mathfrak p$ containing $\overline F$. We get an exact sequence
\[1\to I_{\mathfrak p}\to D_{\mathfrak p}\to\Gal(\overline{\mathbf F}_{\mathfrak p}/\Ffp)
\to1\]
where $D_{\mathfrak p}$ is the decomposition group and
$I_{\mathfrak p}$ the inertia. We denote by
$\widetilde{\Fr}_{\mathfrak p}$ a lifting of the Frobenius map
to $D_{\mathfrak p}\subset\Gal(\overline F/F)$ (which up to conjugation
depends only on $\mathfrak p$), and put
\[L_{\mathfrak p}(s,M)=
\frac{1}{\Det(1-(\cardinal\Ffp)^{-s}\widetilde{\Fr}_{\mathfrak p}\mid
M^{I_{\mathfrak p}})}.\]
We fix a finite set $S$ of bad places containing
the archimedean ones so that $V$ admits a smooth
projective model $\calli V$ over the ring
of $S$-integers $\calli O_S$. For any $\mathfrak p$
in $\Val(F)-S$ we consider
\[L_{\mathfrak p}(s,\Pic(\Vbar))=
L_{\mathfrak p}(s,\Pic(\Vbar)\otimes_\ZZ\QQ).\]
The corresponding global $L$-function is given by the Euler product
\[L_S(s,\Pic(\Vbar))=\prod_{\mathfrak p\in\Val(F)-S}
L_{\mathfrak p}(s,\Pic(\Vbar))\]
which converges for $\re s>1$ and has a meromorphic
continuation to $\CC$ with a pole of order $t=\rk\Pic(V)$
at $1$. One introduces local convergence factors
$\lambda_v$ given by
\[\lambda_v=\begin{cases}
L_v(1,\Pic(\Vbar))\text{ if $v\in\Val(F)-S$},\\
1\text{ otherwise.}
\end{cases}\]
The Weil conjectures (proved by Deligne) imply that the
Tamagawa measure
\[\prod_{v\in\Val(F)}\lambda_v^{-1}\omegaHv\]
converges on $V(\Adeles_F)$ (see \cite[proposition 2.2.2]{peyre:fano}).
\end{notas}

\begin{defi}
The {\em Tamagawa measure} on $V(\Adeles_F)$ corresponding
to the adelic metric $(\metric)_{v\in\Val(F)}$ is defined by
\[\omegaH=\frac{1}{\sqrt {d_F}^{\,\dim V}}\lim_{s\to 1}
(s-1)^tL_S(s,\Pic(\Vbar))\prod_{v\in\Val(F)}\lambda^{-1}_v\omegaHv.\]
\par
\end{defi}

From the arithmetic standpoint, 
it seems more natural to 
integrate $\omegaH$ over
the closure $\overline{V(F)}\subset V(\Adeles_F)$ 
(as in the original approach to the Tamagawa number). 
However, computationally, it is easier to work with 
$V(\Adeles_F)^{\Br}$. Therefore, 
following a suggestion of Salberger, we define
here

\begin{defi}
\begin{equation}
\label{equ:constant:tamagawa}
\tauH(V)=\omegaH(V(\Adeles)^{\Br})
\end{equation}
\end{defi}

Let $\Pic(V)\dual$ be the dual lattice to $\Pic(V)$.
We denote by $\Haar{\boldit y}$ the corresponding Lebes\-gue
measure on $\Pic(V)\dual\otimes_\ZZ\RR$ and by
\[\Ceff(V)\dual=\{\,x\in\Pic(V)\dual\otimes_\ZZ\RR\mid
\forall y\in\Ceff(V),\,\langle x, y\rangle\vargeq 0\,\}\]
the dual cone of $\Ceff(V)$.

\begin{defi}
We define
\[\alpha(V)=\frac{1}{(t-1)!}\int_{\Ceff(V)\dual}
e^{-\langle\antican,\boldit y\rangle}\Haar{\boldit y}\]
and
\[\beta(V)=\cardinal H^1(k,\Pic(\overline{V})).\]
The theoretical constant attached to $V$ and $H$
is defined as
\begin{equation}
\label{equ:constant:constant}
\thetaH(V)=\alpha(V)\beta(V)\tauH(V).
\end{equation}
\end{defi}

In the following sections we compute $\thetaH(V)$ for various
diagonal cubic surfaces.

\section{The Galois module $\Pic(\Vbar)$}
\label{section:lines}

The description of this Galois module is based upon the study
of the 27 lines of the cubic. We fix notations
for these lines extending those given by Colliot-Th\'el\`ene,
Kanevsky and Sansuc in \cite[p. 9]{colliotkanevskysansuc:cubic}.

\begin{notas}
From now on $V$ is a diagonal
cubic surface $V$ given by an equation of the form
\begin{equation}
\label{equ:lines:cubic}
aX^3+bY^3+cZ^3+dT^3=0
\end{equation}
where a,b,c and d are 
strictly positive integers with $\gcd(a,b,c,d)=1$.
Let
\[S=\{\infty,3\}\cup\{\,p\mid\,p\divise abc\}\]
\par
We fix a cubic root $\alpha$ (resp. $\alpha',\alpha''$)
of $b/a$ (resp. $c/a$, $d/a$) (which is assumed to be in
$\QQ$ if $b/a$ (resp. $c/a$, $d/a$) is a cube in $\QQ$)
and we put
\begin{align*}
\beta&=\frac{\alpha''}{\alpha'}=
\sqrt[\uproot{5}\leftroot{-3}3]{\frac{d}{c}},&
\beta'&=\frac{\alpha''}{\alpha}=
\sqrt[\uproot{5}\leftroot{-3}3]{\frac{d}{b}}&&
\text{and}&
\beta''&=\frac{\alpha'}{\alpha}=
\sqrt[\uproot{5}\leftroot{-3}3]{\frac{c}{b}}.\\
\noalign{\noindent We also consider}
\gamma &=\frac{\alpha''}{\alpha\alpha'}
=\sqrt[\uproot{5}\leftroot{-3}3]{\frac{ad}{bc}},&
\gamma' &=\frac{\alpha}{\alpha'\alpha''}
=\sqrt[\uproot{5}\leftroot{-3}3]{\frac{ab}{cd}}&&
\text{and}&
\gamma'' &=\frac{\alpha'}{\alpha\alpha''}
=\sqrt[\uproot{5}\leftroot{-3}3]{\frac{ac}{bd}}.
\end{align*}
We denote by $\theta$ a primitive third root of one.
The 27 lines of the cubic surface \eqref{equ:lines:cubic}
are given by the following equations, where $i$
belongs to $\ZZ/3\ZZ$:
\[
\arraycolsep 0pt
\begin{array}{rclrclrcl}
L(i)&:&\begin{cases}
X{+}\theta^i\alpha Y=0,\\
Z{+}\theta^i\beta T=0.
\end{cases}&
L'(i)&:&\begin{cases}
X{+}\theta^i\alpha Y=0,\\
Z{+}\theta^{i{+}1}\beta T=0.
\end{cases}&
L''(i)&:&\begin{cases}
X{+}\theta^i\alpha Y=0,\\
Z{+}\theta^{i{+}2}\beta T=0.
\end{cases}\\
\noalign{\vskip0.75ex\penalty 1000}
M(i)&:&\begin{cases}
X{+}\theta^i\alpha' Z=0,\\
Y{+}\theta^{i{+}1}\beta' T=0.
\end{cases}&
M'(i)&:&\begin{cases}
X{+}\theta^i\alpha' Z=0,\\
Y{+}\theta^{i{+}2}\beta' T=0.
\end{cases}&
M''(i)&:&\begin{cases}
X{+}\theta^i\alpha' Z=0,\\
Y{+}\theta^{i}\beta' T=0.
\end{cases}\\
\noalign{\vskip0.75ex\penalty 1000}
N(i)&:&\begin{cases}
X{+}\theta^i\alpha'' T=0,\\
Y{+}\theta^{i{+}2}\beta'' Z=0.
\end{cases}&
N'(i)&:&\begin{cases}
X{+}\theta^i\alpha'' T=0,\\
Y{+}\theta^{i}\beta'' Z=0.
\end{cases}&
N''(i)&:&\begin{cases}
X{+}\theta^i\alpha'' T=0,\\
Y{+}\theta^{i{+}1}\beta'' Z=0.
\end{cases}
\end{array}\]
Let $K$ be the field $\QQ(j,\alpha,\alpha',\alpha'').$
It is a Galois extension of $\QQ$. In the generic case, $K$
is an extension of degree 54 with a Galois group isomorphic
to
\[(\ZZ/3\ZZ)^3\rtimes\ZZ/2\ZZ.\]
It is generated by the elements $c$, $\tau$, $\tau'$ and $\tau''$
characterized by their action on
$\theta$, $\alpha$, $\alpha'$ and $\alpha''$.
\[\setbox\strutbox\vbox to 0pt{}
\begin{array}{|c||c|c|c|c|}
\cline{2-5}
\multicolumn{1}{c|}{\vbox{\vphantom{$\Bigl($}}}&
\theta&\alpha&\alpha'&\alpha''\\
\cline{2-5}\multicolumn{1}{c|}{
\vbox{\hrule height\doublerulesep depth 0pt width 0pt}}&&&&\\
\hline
\vbox{\vphantom{$\Bigl($}}c&
\theta^2&\alpha&\alpha'&\alpha''\\
\hline
\vbox{\vphantom{$\Bigl($}}\tau&
\theta&\theta\alpha&\alpha'&\alpha''\\
\hline
\vbox{\vphantom{$\Bigl($}}\tau '&
\theta&\alpha&\theta\alpha'&\alpha''\\
\hline
\vbox{\vphantom{$\Bigl($}}\tau''&
\theta&\alpha&\alpha'&\theta\alpha''\\
\hline
\end{array}\]
\end{notas}
Their action on the 27 lines is given as follows:
for $\tau$ we have
\begin{align}
\label{equ:lines:tau}
&\begin{array}{rcl}
L(i)&\mathclap{\longrightarrow}&L''(i+1)\\
\nwarrow&&\swarrow\\
&\mathclap{L'(i+2)}
\end{array}&&
\begin{array}{rcl}
M(i)&\mathclap{\longrightarrow}&M''(i)\\
\nwarrow&&\swarrow\\
&\mathclap{M'(i)}
\end{array}&&
\text{and}&&
\begin{array}{rcl}
N(i)&\mathclap{\longrightarrow}&N''(i)\\
\nwarrow&&\swarrow\\
&\mathclap{N'(i)}
\end{array}
\\
\noalign{\noindent for $\tau'$:}
\label{equ:lines:tauprime}
&\begin{array}{rcl}
L(i)&\mathclap{\longrightarrow}&L''(i)\\
\nwarrow&&\swarrow\\
&\mathclap{L'(i)}
\end{array}&&
\begin{array}{rcl}
M(i)&\mathclap{\longrightarrow}&M''(i+1)\\
\nwarrow&&\swarrow\\
&\mathclap{M'(i+2)}
\end{array}&&
\text{and}&&
\begin{array}{rcl}
N(i)&\mathclap{\longrightarrow}&N'(i)\\
\nwarrow&&\swarrow\\
&\mathclap{N''(i)}
\end{array}
\\
\noalign{\noindent for $\tau''$:}
\label{equ:lines:tausecond}
&\begin{array}{rcl}
L(i)&\mathclap{\longrightarrow}&L'(i)\\
\nwarrow&&\swarrow\\
&\mathclap{L''(i)}
\end{array}&&
\begin{array}{rcl}
M(i)&\mathclap{\longrightarrow}&M'(i)\\
\nwarrow&&\swarrow\\
&\mathclap{M''(i)}
\end{array}&&
\text{and}&&
\begin{array}{rcl}
N(i)&\mathclap{\longrightarrow}&N''(i+1)\\
\nwarrow&&\swarrow\\
&\mathclap{N'(i+2)}
\end{array}\end{align}
for $c$:
\begin{equation}
\label{equ:lines:c}
\def\p(#1,#2)#3:#4{\rput(#1,#2){\rnode
{#3}{\psframebox*[linestyle=none]{$#4$}}}}
\def\l(#1,#2){\ncline[linewidth=0.01]{<->}{#1}{#2}}
\psset{xunit=1.50cm}
\begin{gathered}
\begin{pspicture}(-.2,-.5)(6,3)
\p(.5,2.5)AA:{L(0)}\p(1.5,2.5)BA:{L'(0)}\p(2.5,2.5)CA:{L''(0)}
\p(.5,1.5)AB:{L(1)}\p(1.5,1.5)BB:{L'(1)}\p(2.5,1.5)CB:{L''(1)}
\p(.5,.5)AC:{L(2)}\p(1.5,.5)BC:{L'(2)}\p(2.5,.5)CC:{L''(2)}
\p(3.5,2.5)DA:{M(0)}\p(4.5,2.5)EA:{M'(0)}\p(5.5,2.5)FA:{M''(0)}
\p(3.5,1.5)DB:{M(1)}\p(4.5,1.5)EB:{M'(1)}\p(5.5,1.5)FB:{M''(1)}
\p(3.5,.5)DC:{M(2)}\p(4.5,.5)EC:{M'(2)}\p(5.5,.5)FC:{M''(2)}
\l(BA,CA)\l(AB,AC)\l(BB,CC)\l(BC,CB)
\l(DA,EA)\l(FB,FC)\l(DB,EC)\l(DC,EB)
\end{pspicture}\\
\begin{pspicture}(-.2,-.5)(3,3)
\p(.5,2.5)AA:{N(0)}\p(1.5,2.5)BA:{N'(0)}\p(2.5,2.5)CA:{N''(0)}
\p(.5,1.5)AB:{N(1)}\p(1.5,1.5)BB:{N'(1)}\p(2.5,1.5)CB:{N''(1)}
\p(.5,.5)AC:{N(2)}\p(1.5,.5)BC:{N'(2)}\p(2.5,.5)CC:{N''(2).}
\nccurve[angleA=15,angleB=165,linewidth=0.01]{<->}{AA}{CA}
\l(BB,BC)\l(AB,CC)\l(AC,CB)
\end{pspicture}
\end{gathered}
\end{equation}
To describe the relations between the classes of these divisors
in $\Pic(\Vbar)$ we consider $\Vbar$ as the blow-up of a plane
$\mathbf P^2_{\overline{\QQ}}$ in six points
$P_1$, $P_2$, $P_3$, $P_4$, $P_5$ and $P_6$. The
27 lines may be described as the $6$ exceptional divisors
$E_1,\dots,E_6$, the 15 strict transforms $L_{i,j}$ of the
projective lines $(P_iP_j)$ for $1\varleq i<j\varleq 6$ and
the 6 strict transforms of the conics $Q_i$ going through all
points except $P_i$. Let $\Lambda$ be the preimage of a line
of $\mathbf P^2_{\overline{\QQ}}$ which does not contain any of the
points $P_1,\dots,P_6$. Then
\[([\Lambda],[E_1],[E_2],[E_3],[E_4],[E_5],[E_6])\]
is a basis of $\Pic(\Vbar)$ and we have the following relations in
$\Pic(\Vbar)$:
\begin{equation}
\label{equ:lines:relations}
\begin{aligned}{}
[L_{i,j}]&=[\Lambda]-[E_i]-[E_j]\qquad\text{for $1\varleq i<j\varleq 6$},\\
[Q_i]&=2[\Lambda]-\sum_{j\neq i}[E_j].
\end{aligned}
\end{equation}
In the following, we choose the projection of $\Vbar$ to
$\mathbf P^2_{\overline{\QQ}}$ so that we have the equalities:
\begin{equation}
\label{equ:lines:bijection}
\begin{aligned}
E_1&=L(0),&E_2&=L(1),&E_3&=L(2),\\
E_4&=M(0),&E_5&=M(1),&E_6&=M(2),\\
Q_1&=L'(1),&Q_2&=L'(2),&Q_3&=L'(0),\\
Q_4&=M''(2),\quad&Q_5&=M''(0),\quad&Q_6&=M''(1),\\
L_{1,2}&=L''(1),&L_{2,3}&=L''(2),&L_{3,1}&=L''(0),\\
L_{4,5}&=M'(0),&L_{5,6}&=M'(1),&L_{6,4}&=M'(2),\\
L_{1,4}&=N'(2),&L_{1,5}&=N'(0),&L_{1,6}&=N'(1),\\
L_{2,4}&=N''(0),&L_{2,5}&=N''(1),&L_{2,6}&=N''(2),\\
L_{3,4}&=N(1),&L_{3,5}&=N(2),&L_{3,6}&=N(0).
\end{aligned}
\end{equation}
\begin{notas}
\label{notas:lines:etale}
We consider the \'etale algebra $E_1$ over $\QQ$
defined as $\QQ(\gamma)$ if $ad/bc$ is not a cube in $\QQ$ and
as $\QQ(\theta)\times\QQ$ otherwise. Similarly, we define
the algebra $E_2$ (resp. $E_3$) corresponding to
$\gamma'$ (resp. $\gamma''$) and we put
\[E=E_1\times E_2\times E_3.\]
We also consider the following elements of $\Pic(\Vbar)$
\begin{align*}
e^1_0&=[L(0)]+[L(1)]+[L(2)],&e^1_1&=[L'(0)]+[L'(1)]+[L'(2)],\\
e^1_2&=[L''(0)]+[L''(1)]+[L''2],
&e^2_0&=[M''(0)]+[M(1)]+[M'(2)],\\
e^2_1&=[M'(0)]+[M''(1)]+[M(2)],
&e^2_2&=[M(0)]+[M'(1)]+[M''(2)],\\
e^3_0&=[N'(0)]+[N'(1)]+[N'(2)],&e^3_1&=[N''(0)]+[N''(1)]+[N''(2)],\\
e^3_2&=[N(0)]+[N(1)]+[N(2)]
\end{align*}
and the sets
\[\calli E_1=\{e^1_0,e^1_1,e^1_2\},\quad
\calli E_2=\{e^2_0,e^2_1,e^2_2\},\quad
\calli E_3=\{e^3_0,e^3_1,e^3_2\},\]
and $\calli E=\calli E_1\sqcup\calli E_2\sqcup\calli E_3$.
\end{notas}
\begin{lem}
The sets $\calli E_1$, $\calli E_2$ and $\calli E_3$
are globally invariant under the action of $\Gal(K/\QQ)$
and the \'etale algebra corresponding to the set $\calli E_i$
is isomorphic to $E_i$.
\end{lem}
\begin{proof}
The fact that the sets  $\calli E_1$, $\calli E_2$ and $\calli E_3$
are globally invariant follows immediately
from the descriptions \eqref{equ:lines:tau}--\eqref{equ:lines:c}.
The \'etale algebra $F$ corresponding to a 
finite $\Gal(K/\QQ)$-set $\calli F$
may be defined as the algebra
\[(K[\calli F])^{\Gal(K/\QQ)}\]
where $K[\calli F]$ is the algebra $K^{\calli F}$ and
where $\Gal(K/\QQ)$ acts simultaneously on $K$ and $\calli F$.
In the generic case, let us consider
\[\sigma=\tau\tau',\quad\sigma'=\tau'\tau''\quad
\text{and}\quad\sigma''=\tau''\tau.\]
Then $\sigma$ sends $\gamma$ on $\theta\gamma$ and acts trivially
on $\gamma'$, $\gamma''$ and $\theta$. We may describe similarly
the actions of $\sigma'$ and $\sigma''$.
The action of $\Gal(K/\QQ)$ on $\calli E_1$ in the generic case
is given by the table
\[\setbox\strutbox\vbox to 0pt{}
\begin{array}{|c||c|c|c|}
\cline{2-4}
\multicolumn{1}{c|}{\vbox{\vphantom{$\Bigl($}}}&
e^1_0&e^1_1&e^1_2\\
\cline{2-4}\multicolumn{1}{c|}{
\vbox{\hrule height\doublerulesep depth 0pt width 0pt}}&&&\\
\hline
\vbox{\vphantom{$\Bigl($}}c&
e^1_0&e^1_2&e^1_1\\
\hline
\vbox{\vphantom{$\Bigl($}}\sigma&
e^1_1&e^1_2&e^1_0\\
\hline
\vbox{\vphantom{$\Bigl($}}\sigma '&
e^1_0&e^1_1&e^1_2\\
\hline
\vbox{\vphantom{$\Bigl($}}\sigma''&
e^1_0&e^1_1&e^1_2\\
\hline
\end{array}\]
This implies that if $ad/bc$ is not a cube in $\QQ$, then
$\calli E_1$ is isomorphic to
\[\Gal(K/\QQ)/\Gal(K/\QQ(\gamma))\]
as a $\Gal(K/\QQ)$-set. Then the corresponding \'etale algebra
is
\[(K[\Gal(K/\QQ)/\Gal(K/\QQ(\gamma))])^{\Gal(K/\QQ)}\iso
K^{\Gal(K/\QQ(\gamma))}=\QQ(\gamma)=E_1.\]
Similarly if $ad/bc$ is a cube in $\QQ$, then we may decompose
$\calli E_1$ into two orbits and we see that the corresponding \'etale
algebra is
$\QQ(\theta)\times\QQ=E_1.$
The proofs for $\calli E_2$ and $\calli E_3$ are similar.
\end{proof}
\begin{lem}
\label{lem:lines:exact}
There exists an exact sequence of $\Gal(K/\QQ)$ modules
\[0\to\QQ^2\to\QQ[\calli E]\to\Pic(\Vbar)\otimes_\ZZ\QQ\to 0.\]
\end{lem}
\begin{proof}
By \eqref{equ:lines:relations} and \eqref{equ:lines:bijection},
we have in $\Pic(\Vbar)$ the relations
\begin{align*}
e^1_0&=[E_1]+[E_2]+[E_3],\\
e^1_1&=6[\Lambda]-2[E_1]-2[E_2]-2[E_3]-3[E_4]
-3[E_5]-3[E_6],\\
e^1_2&=3[\Lambda]-2[E_1]-2[E_2]-2[E_3],\\
e^2_0&=3[\Lambda]-[E_1]-[E_2]-[E_3]-2[E_4]
+[E_5]-2[E_6],\\
e^2_1&=3[\Lambda]-[E_1]-[E_2]-[E_3]-2[E_4]
-2[E_5]+[E_6],\\
e^2_2&=3[\Lambda]-[E_1]-[E_2]-[E_3]+[E_4]
-2[E_5]-2[E_6],\\
e^3_0&=3[\Lambda]-3[E_1]-[E_4]-[E_5]-[E_6],\\
e^3_1&=3[\Lambda]-3[E_2]-[E_4]-[E_5]-[E_6],\\
e^3_2&=3[\Lambda]-3[E_3]-[E_4]-[E_5]-[E_6]
\end{align*}
which proves that the natural projection from
$\QQ[\calli E]$ to $\Pic(\Vbar)\otimes_\ZZ\QQ$ is surjective. Moreover
one has the relations
\[3\antican=\sum_{x\in\calli E_1}x=
\sum_{x\in\calli E_2}x=\sum_{x\in\calli E_3}x,\]
which gives a homomorphism of $\Gal(K/\QQ)$-modules
\[\QQ^2\to\QQ[\calli E]\]
and the exact sequence of the lemma.
\end{proof}
\begin{notas}
For any prime $p$ and any finite field extension $F$ of $\QQ$,
we consider the local factor $\zeta_{F,p}$
of the function $\zeta_F$ at $p$
which is defined by
\[\zeta_{F,p}(s)=\prod_{\{\,v\in\Val(F)\mid v\divise p\,\}}
(1-\cardinal \mathbf F^{-s}_v)^{-1}.\]
Let $F$ be an \'etale algebra over $\QQ$ and $F=\prod_{i\in I}F_i$
its decomposition in fields. Put
\[\zeta_F(s)=\prod_{i\in I}\zeta_{F_i}(s)\quad
\text{and}\quad
\zeta_{F,p}(s)=\prod_{i\in I}\zeta_{F_i,p}(s).\]
For any prime $p$, we denote by $\nu_F(p)$ the number of components
of $F\otimes_\QQ\QQ_p$ of degree one over $\QQ_p$.
\end{notas}
\begin{prop}
\label{prop:lines:l}
With notation as above,
for any prime $p$ not in $S$, one has
\begin{itemize}
\item[(i)]$\displaystyle L_p(s,\Pic(\Vbar))=
\frac{\zeta_{E,p}(s)}{\zeta_{\QQ,p}(s)^2},$
\item[(ii)]$\displaystyle\Trace(\widetilde\Fr_p\mid\Pic(\Vbar))=\nu_E(p)-2$.
\end{itemize}
\end{prop}
\begin{proof}
By lemma \ref{lem:lines:exact}, we have
\[L_p(s,\Pic(\Vbar))=\frac{L_p(s,\QQ[\calli E])}{L_p(s,\QQ)^2}.\]
Thus it is enough to prove that if $E$ is an arbitrary \'etale algebra
over $\QQ$ corresponding to a $\Gal(\overline{\QQ}/\QQ)$-set
$\calli E$ and if $p$ is a prime such that $E/\QQ$ is not ramified
at $p$, then
\[\zeta_{E,p}(s)=L_p(s,\QQ[\calli E]).\]
This well-known assertion follows from the fact
that the components of $E\otimes\QQ_p$ are in bijection
with the orbits of $\Fr_p$ in $\calli E$,
and the degree of each component is the length 
of the corresponding orbit. This proves (i).
\par
But this also shows that
\[\Trace(\widetilde\Fr_p\mid\QQ[\calli E])
=\nu_E(p)\]
which implies (ii).
\end{proof}
\begin{rem}
Thus the factor $\lambda'_p$ which was defined
in proposition 5.1 in \cite{peyretschinkel:numericone}
coincides with $L_p(1,\Pic(\Vbar))$ at the good places
(as suggested by the referee of that paper).
\end{rem}

\section{Euler product for the good places}
\label{section:euler}

We need to compute the number of solutions of \eqref{equ:lines:cubic}
modulo $p$ for all primes not in $S$.

\begin{prop}
For any prime $p$ not in $S$, one has
\[\frac{\cardinal V(\Fp)}{p^2}=1+\frac{\nu_E(p)-2}{p}+\frac{1}{p^2}\]
where $E$ is the \'etale algebra defined in \S\ref{section:lines}.
\end{prop}
\begin{proof}
By a result of Weil (see \cite[theorem 23.1]{manin:cubic}), 
\[\cardinal V(\Fp)=1+\Trace(\Fr_p\mid\Pic(\Vbar))p+p^2.\]
Proposition \ref{prop:lines:l} implies that
\[\Trace(\Fr_p\mid\Pic(\Vbar))=\nu_E(p)-2.\qed\]
\noqed
\end{proof}
\begin{rem}
We could have proved this result directly as in 
\cite{peyretschinkel:numericone}. Let $N(p)$
be the number of solutions of \eqref{equ:lines:cubic}
in $\Fp^4$. We have
\[\cardinal V(\Fp)=\frac{N(p)-1}{p-1}.\]
By 
\cite[\S8.7 theorem 5]{irelandrosen:number}, one has
\[N(p)=p^3+\sum\overline\chi_1(a)\overline\chi_2(b)
\overline\chi_3(c)\overline\chi_4(d)
J_0(\chi_1,\chi_2,\chi_3,\chi_4),\]
where the sum is taken over all  quadruples $(\chi_1,\dots,\chi_4)$
of nontrivial cubic characters from $\Fp^*$ to $\CC^*$ such that
$\chi_1\chi_2\chi_3\chi_4=1$ and where
\[J_0(\chi_1,\chi_2,\chi_3,\chi_4)=\sum_{t_1+\dots+t_4=0}
\prod_{i=1}^4\chi_i(t_i),\]
the characters being extended by $\chi_i(0)=0$. For $p\equiv 2$ mod $3$
there are no nontrivial characters and the formula is obvious. Otherwise
there are exactly two nontrivial conjugated characters $\chi$ and $\overline\chi$. 
By \cite[proof of prop. 4.1]{peyretschinkel:numericone},
we have
\[J_0(\chi_1,\chi_2,\chi_3,\chi_4)=p(p-1)\]
and
\[\cardinal V(\Fp)=1+p(1+\sum\chi_1(a)\chi_2(b)\chi_3(c)\chi_4(d))+p^2\]
where the sum is taken over the same quadruples as above. The formula
\begin{multline*}
\sum\chi_1(a)\chi_2(b)\chi_3(c)\chi_4(d)=\\
\chi\Bigl(\frac{ab}{cd}\Bigr)+\overline\chi\Bigl(\frac{ab}{cd}\Bigr)+
\chi\Bigl(\frac{ac}{bd}\Bigr)+\overline\chi\Bigl(\frac{ac}{bd}\Bigr)+
\chi\Bigl(\frac{ad}{bc}\Bigr)+\overline\chi\Bigl(\frac{ad}{bc}\Bigr)
\end{multline*}
implies the result.
\end{rem}
\begin{notas}
For any place $v$ of $\QQ$, we put
\[\lambda_v=\begin{cases}
\frac{\zeta_{E,v}(s)}{\zeta_{\QQ,v}(s)^2}\text{ if $v$ is finite},\\
1\text{ otherwise.}
\end{cases}\]
\end{notas}
\begin{rem}
\label{rem:euler:euler}
By proposition \ref{prop:lines:l},
$\lambda_p=L_p(1,\Pic(\Vbar))$ if $p\in\Val(\QQ)-S$.
Thus the Tamagawa measure $\omegaH$ is given by the formula
\[\omegaH=\lim_{s\to 1}(s-1)^{\rk\Pic(V)}\Bigl(
\frac{\zeta_E(s)}{\zeta_\QQ(s)^2}\Bigr)\times
\prod_{v\in\Val(\QQ)}\lambda_v^{-1}\omegaHv.\]
By lemmata 3.2 and 3.4 in \cite{peyretschinkel:numericone}
and lemma 5.4.6 in \cite{peyre:fano}, for any $p$
in $\Val(\QQ)-S$ one has
\[\omegaHp(V(\QQ_p))=\frac{\cardinal V(\Fp)}{p^2}\]
(see also \cite[lemma 2.2.1]{peyre:fano} and
\cite[remark 5.2]{peyretschinkel:numericone}).
Therefore, the local factor at a good place $p$ is given by
\[\begin{array}{cl}
\Bigl(1-\frac 1p\Bigr)^7\Bigl(1+\frac 7p+\frac 1{p^2}\Bigr)&
\text{if $p\equiv 1$ mod $3$ and $\nu_E(p)=9$}\\
\Bigl(1-\frac 1p\Bigr)^4\Bigl(1-\frac 1{p^3}\Bigr)
\Bigl(1+\frac 4p+\frac 1{p^2}\Bigr)&
\text{if $p\equiv 1$ mod $3$ and $\nu_E(p)=6$}\\
\Bigl(1-\frac 1p\Bigr)\Bigl(1-\frac 1{p^3}\Bigr)^2
\Bigl(1+\frac 1p+\frac 1{p^2}\Bigr)&
\text{if $p\equiv 1$ mod $3$ and $\nu_E(p)=3$}\\
\Bigl(1-\frac 1p\Bigr)^{-2}\Bigl(1-\frac 1{p^3}\Bigr)^3
\Bigl(1-\frac 2p+\frac 1{p^2}\Bigr)&
\text{if $p\equiv 1$ mod $3$ and $\nu_E(p)=0$}\\
\Bigl(1-\frac 1p\Bigr)\Bigl(1-\frac 1{p^2}\Bigr)^3
\Bigl(1+\frac 1p+\frac 1{p^2}\Bigr)&
\text{if $p\equiv 2$ mod $3$.}
\end{array}\]
We get (for the good places) the factors $C_0$, $C_1$, $C_2$
and $C_3$ where
\begin{align*}
C_0&=\prod_{
\begin{scriptarray}
p\notdivise 3abcd,\crr
p\equiv 2\text{ mod }3.
\end{scriptarray}}
\Bigl(1-\frac 1{p^3}\Bigr)\Bigl(1-\frac 1{p^2}\Bigr)^3,\\
C_1&=\prod_{
\begin{scriptarray}
p\notdivise 3abcd,\crr
p\equiv 1\text{ mod }3,\crr
\nu_E(p)=9.
\end{scriptarray}}
\Bigl(1-\frac 1p\Bigr)^7\Bigl(1+\frac 7p+\frac 1{p^2}\Bigr),\\
C_2&=\prod_{
\begin{scriptarray}
p\notdivise 3abcd,\crr
p\equiv 1\text{ mod }3,\crr
\nu_E(p)=6.
\end{scriptarray}}
\Bigl(1-\frac 1{p^3}\Bigr)
\Bigl(1-\frac 1p\Bigr)^4\Bigl(1+\frac 4p+\frac 1{p^2}\Bigr),\\
C_3&=\prod_{
\begin{scriptarray}
p\notdivise 3abcd,\crr
p\equiv 1\text{ mod }3,\crr
\nu_E(p)=0\text{ or }3.
\end{scriptarray}}
\Bigl(1-\frac 1{p^3}\Bigr)^3.
\end{align*}
These products converge rapidly and are easily approximated.
\end{rem}

\section{Density at the bad places}
\label{section:bad}
In this section we restrict to cubic surfaces with equations
of the form
\begin{equation}
\label{equ:bad:rktwo}
X^3+Y^3+qZ^3+q^2T^3=0
\end{equation}
with $q$ prime and
\begin{equation}
\label{equ:bad:rkthree}
aX^3+aY^3+qZ^3+qT^3=0
\end{equation}
with $q$ prime and $a$ an integer coprime to $q$.
\begin{notas}
If $V$ is defined by the equation \eqref{equ:lines:cubic},
and $p$ is a prime, then we consider
\[N^*(p^r)=\cardinal\{(x,y,z,t)\in(\ZZ/p^r\ZZ)^4
{-}(p\ZZ/p^r\ZZ)^4\mid
ax^3{+}by^3{+}cz^3{+}dt^3{=}0\text{ in }\ZZ/p^r\ZZ\}\]
\end{notas}
\begin{rem}
By \cite[lemmata 3.2 and 3.4]{peyretschinkel:numericone},
there is an explicit integer $r_0$ such that
\[\omegaHp(V(\QQ_p))=\frac1{1-p^{-1}}\times
\frac{N^*(p^{r_0})}{p^{3r_0}}.\]
If $p=3$ and $3\notdivise abcd$, then a direct computation
in $(\ZZ/9\ZZ)^4$ gives the value of $N^*(9)$ and
thus of $\omegaHp(V(\QQ_p))$.
Thus, in the following lemma we restrict 
to the case when $V$ is given by \eqref{equ:bad:rktwo}
or \eqref{equ:bad:rkthree} and $p=q$.
\end{rem}
\begin{lem}
If $V$ is given by the equation
\[X^3+Y^3+pZ^3+p^2T^3=0\]
then for $r\vargeq 2$,
\[\frac{N^*(p^r)}{p^{3r}}=\begin{cases}
1-\frac 1p&\text{ if $p\equiv 2$ mod $3$},\\
3\Bigl(1-\frac 1p\Bigr)&\text{ if $p\equiv 1$ mod $3$},\\
\frac 23&\text{ if $p=3$}.
\end{cases}\]
If $V$ is given by the equation
\[aX^3+aY^3+pZ^3+pT^3=0,\]
with $p\notdivise a$, then for $r\vargeq 3$,
\[\frac{N^*(p^r)}{p^{3r}}=\begin{cases}
1-\frac 1{p^2}&\text{ if $p\equiv 2$ mod $3$},\\
3\Bigl(1-\frac 1{p^2}\Bigr)&\text{ if $p\equiv 1$ mod $3$},\\
\frac 43&\text{ if $p=3$}.
\end{cases}\]
\end{lem}
\begin{rem}
This lemma implies that if $V$ is given by the 
first equation then the local factor at $p$ is given by
\[\lambda_p\omegaHp(V(\QQ_p))=
\begin{cases}
\Bigl(1-\frac 1{p^2}\Bigr)\Bigl(1-\frac 1p\Bigr)&
\text{if $p\equiv 2$ mod $3$},\\
3\Bigl(1-\frac 1p\Bigr)^3&
\text{if $p\equiv 1$ mod $3$},\\
\frac 49&
\text{if $p=3$},
\end{cases}\]
and if $V$ is given by the second equation
then this factor is
\[\lambda_p\omegaHp(V(\QQ_p))=\begin{cases}
\Bigl(1-\frac 1{p^2}\Bigr)^3&\text{if $p\equiv 2$ mod 3},\\
3\Bigl(1-\frac 1p\Bigr)^4\Bigl(1-\frac 1{p^2}\Bigr)&
\text{if $p\equiv 1$ mod $3$},\\
\frac {16}{27}&\text{if $p=3$}.
\end{cases}\]
\end{rem}
\begin{proof}
Let us consider the set of quadruples $(x,y,z,t)$
in $(\ZZ/p^r\ZZ)^4-(p\ZZ/p^r\ZZ)^4$ such that
\begin{equation}
\label{equ:bad:cubicone}
x^3+y^3+pz^3+p^2t^3=0\qquad\text{in $\ZZ/p^r\ZZ$.}
\end{equation}
If $p\divise x$ then $p\divise y$, $p\divise z$ and $p\divise t$.
Therefore, for any $(x,y,z,t)$ as above, $p\notdivise x$ and
$p\notdivise y$. But for any triple
$(y,z,t)$ in $(\ZZ/p^r\ZZ-p\ZZ/p^r\ZZ)\times(\ZZ/p^r\ZZ)^2$, there
exists exactly one $x$ verifying \eqref{equ:bad:cubicone} if
$p\equiv 2$ mod $3$ and exactly three of them if $p\equiv 1$ mod $3$.
If $p=3$ and $y$ belongs to $\ZZ/3^r\ZZ-3\ZZ/3^r\ZZ$
then \eqref{equ:bad:cubicone} implies that $3\divise z$.
For any triple $(y,z,t)$ with $y$ in $(\ZZ/3^r\ZZ)-(3\ZZ/3^r\ZZ)$,
$z$ in $(3\ZZ/3^r\ZZ)$ and $t$
in $(\ZZ/3^r\ZZ)$ there exists exactly three $x$ in $\ZZ/3^r\ZZ$
which satisfy \eqref{equ:bad:cubicone}. We get that
\[\frac{N^*(p^r)}{p^{3r}}=\begin{cases}
\frac{(p-1)p^{r-1}\times p^r\times p^r}{p^{3r}}=1-\frac 1p
&\text{if $p\equiv 2$ mod $3$},\\
3\frac{(p-1)p^{r-1}\times p^r\times p^r}{p^{3r}}=3\Bigl(1-\frac 1p\Bigr)
&\text{if $p\equiv 1$ mod $3$},\\
3\frac{2\times 3^{r-1}\times 3^{r-1}\times 3^r}{3^{3r}}=\frac 23&
\text{if $p=3$.}
\end{cases}\]
Let us now turn to the set of $(x,y,z,t)$ in $(\ZZ/p^r\ZZ)^4-
(p\ZZ/p^r\ZZ)^4$ such that
\[ax^3+ay^3+pz^3+pt^3=0.\]
We decompose this set as follows
\begin{align*}
N_1^*(p^r)&=\cardinal\left\{\,(x,y,z,t)\in(\ZZ/p^r\ZZ)^4
{-}(p\ZZ/p^r\ZZ)^4\left\vert\begin{cases}
p\notdivise x,\\
ax^3+ay^3+pz^3+pt^3=0.
\end{cases}
\,\right.\right\}\\
N_2^*(p^r)&=\cardinal\left\{\,(x,y,z,t)\in(\ZZ/p^r\ZZ)^4
{-}(p\ZZ/p^r\ZZ)^4\left\vert\begin{cases}
p\divise x,\quad p\notdivise z,\\
ax^3+ay^3+pz^3+pt^3=0.
\end{cases}
\,\right.\right\}
\end{align*}
As above we have the formula
\[\frac{N^*_1(p^r)}{p^{3r}}=
\begin{cases}
\frac{(p-1)p^{r-1}\times p^r\times p^r}{p^{3r}}=1-\frac 1p
&\text{if $p\equiv 2$ mod $3$,}\\
3\times \frac{(p-1)p^{r-1}\times p^r\times p^r}{p^{3r}}=
3\Bigl(1-\frac 1p\Bigr)
&\text{if $p\equiv 1$ mod $3$,}\\
3\frac{2\times 3^{r-1}\times 3^{r-1}\times 3^r}{3^{3r}}=\frac 23&
\text{if $p=3$,}
\end{cases}\]
where for $p=3$ we use the equality
\[3^{r-1}\times 3^r=\cardinal\{\,(z,t)\in(\ZZ/3^r\ZZ)^2\mid
z^3\equiv t^3\text{ mod }3\,\}.\]
On the other hand,
\[N_2^*(p^r)=p^2\left\{\,(x,y,z,t)\in(\ZZ/p^{r-1}\ZZ)^4\left\vert
\begin{cases}
p\notdivise z\\
ap^2x^3+ap^2y^3+z^3+t^3=0.
\end{cases}
\,\right.\right\}\]
and
\[\frac{N^*_2(p^r)}{p^{3r}}=
\frac{p^2}{p^3}\times\begin{cases}
\frac{(p-1)p^{r-2}\times p^{r-1}\times p^{r-1}}{p^{3(r-1)}}=1-\frac 1p
&\text{if $p\equiv 2$ mod $3$},\\
3\frac{(p-1)p^{r-2}\times p^{r-1}\times p^{r-1}}{p^{3(r-1)}}=
3\Bigl(1-\frac 1p\Bigr)
&\text{if $p\equiv 2$ mod $3$},\\
3\frac{2\times 3^{r-2}\times 3^{r-1}\times 3^{r-1}}{3^{3(r-1)}}=2
&\text{if $p=3$}.
\end{cases}\]
We conclude:
\[\frac{N^*(p^r)}{p^{3r}}=
\begin{cases}
1-\frac 1p+\frac 1p-\frac 1{p^2}=1-\frac 1{p^2}&\text{if $p\equiv 2$ mod 3},\\
3\Bigl(1-\frac{1}{p^2}\Bigr)&\text{if $p\equiv 1$ mod 3},\\
\frac 23+\frac 23=\frac43&\text{if $p=3$.}\qed
\end{cases}\]
\noqed
\end{proof}

\section{The constant $\alpha(V)$}
\label{section:alpha}
Since the cubic surfaces we consider in this paper are 
$\QQ$-rational (which implies that $\beta(V)=1$), it
remains to compute the rank $t$ of the Picard group and
the value of $\alpha(V)$.

\begin{prop}
If $V$ is given by the equation
\begin{equation}
\label{equ:alpha:rktwo}
X^3+Y^3+aZ^3+aT^3=0,
\end{equation}
where $a$ is not a cube in $\QQ$, then $\rk\Pic(V)=2$
and $\alpha(V)=2$.
\par
If $V$ is given by the equation
\begin{equation}
\label{equ:alpha:rkthree}
aX^3+aY^3+bZ^3+bT^3=0,
\end{equation}
where $a$ and $b$ are strictly positive integers
and $b/a$ is not a cube in $\QQ$, then $\rk\Pic(V)=3$
and $\alpha(V)=1$.
\par
If $V$ is given by the equation
\begin{equation}
\label{equ:alpha:rkfour}
X^3+Y^3+Z^3+T^3=0
\end{equation}
then $\rk\Pic(V)=4$ and $\alpha(V)=7/18$.
\end{prop}
\begin{proof}
To compute $\alpha(V)$ we shall use its original definition
\cite[\S2]{peyre:fano}: 
\[\alpha(V)=\Vol\{\,x\in\Ceff(V)\mid\langle\antican,x\rangle=1\,\}\]
where the Lebesgue measure on the affine hyperplane
\[\calli H(\lambda)=\{\,x\in\Pic(V)\dual\otimes_\ZZ\RR
\mid\langle\antican,x\rangle
=\lambda\,\}\]
is defined by the $(t-1)$-form $\Haar{\boldit x}$ 
which is characterized by the relation
\[\Haar{\boldit x}\wedge\Haar{\antican}=\Haar{\boldit y}\]
(where $\Haar{\antican}$ is the linear form defined by $\antican$
on $\Pic(V)\dual$ and $\Haar{\boldit y}$ is the form corresponding
to the natural Lebesgue measure on $\Pic(V)\dual\otimes_\ZZ\RR$).
More explicitely, let $(e_1,\dots,e_t)$ be a basis of $\Pic(V)$
and $(e_1\dual,\dots,e_t\dual)$ be the dual basis. Write
\[\antican=\sum_{i=1}^t\lambda_ie_i\]
with $\lambda_t\neq 0$. Let $f_1,\dots,f_{t-1}$ 
be the projection of
$e_1\dual ,\dots,e_{t-1}\dual$ on $\calli H(0)$ along $e_t\dual$. Then
\[\Haar{\boldit x}=\frac 1{\lambda_t}\Haar {f_1\dual}\wedge\dots\wedge
\Haar{f_{t-1}\dual}.\]
\par
When $V$ is given by the equation \eqref{equ:alpha:rktwo} the Galois group
$\Gal(K/\QQ)$ is
\[\ZZ/3\ZZ\rtimes\ZZ/2\ZZ\]
and the orbits of its action on the 27 lines are
\begin{align*}
O_1&=\{L(0),L'(0),L''(0)\},\\
O_2&=\{L(1),L(2),L'(1),L'(2),L''(1),L''(2)\},\\
O_3&=\{M(0),M(2),M'(0),M'(1),M''(1),M''(2)\},\\
O_4&=\{M(1),M'(2),M''(0)\},\\
O_5&=\{N(0),N(1),N'(1),N'(2),N''(0),N''(2)\},\\
O_6&=\{N(2),N'(0),N''(1)\}.
\end{align*}
In the basis $([\Lambda],[E_1],\dots,[E_6])$, a basis of 
$\Pic(V)=(\Pic\Vbar)^{\Gal(K:\QQ)}$
is given by
\[
e_1=\antican,\qquad
e_2=[E_4]-2[E_5]+[E_6].
\]
The effective cone $\Ceff(V)$ is generated by the classes
$[O_i]=\sum_{x\in O_i}x$, which in the basis $(e_0,e_1)$ are given by
\begin{align*}
[O_1]&=e_1,&[O_2]&=2e_1,&[O_3]&=2e_1+e_2,\\
[O_4]&=e_1-e_2,&[O_5]&=2e_1-e_2,&[O_6]&=e_1+e_2.
\end{align*}
Therefore, this cone is generated by the elements $e_1-e_2$ and
$e_1+e_2$ and $\alpha(V)$ is given as the volume of the domain
\[x=1,\quad x+y>0\quad\text{and}\quad x-y>0,\]
that is, as the volume of the segment $[-1,1]$ and
$\alpha(V)=2$.
\par
If $V$ is given by the equation \eqref{equ:alpha:rkthree}
then $\Gal(K/\QQ)$ is isomorphic to
\[\ZZ/3\ZZ\rtimes\ZZ/2\ZZ\]
and the orbits of the Galois action on the $27$ lines are
\begin{align*}
O_1&=\{L(0)\},\\
O_2&=\{L(1),L(2)\},\\
O_3&=\{L'(0),L''(0)\},\\
O_4&=\{L'(1),L''(2)\},\\
O_5&=\{L'(2),L''(1)\},\\
O_6&=\{M(0),M(1),M(2),M'(0),M'(1),M'(2)\},\\
O_7&=\{M''(0),M''(1),M''(2)\},\\
O_8&=\{N(0),N(1),N(2),N''(0),N''(1),N''(2)\},\\
O_9&=\{N'(0),N'(1),N'(2)\}.
\end{align*}
A basis of $\Pic(V)$ is given by
\[e_1=\antican,\quad e_2=[E_1],\quad e_3=[E_2]+[E_3],\]
and the cone $\Ceff(V)$ is generated by
\begin{align*}
[O_1]&=e_2,&[O_2]&=e_3,&[O_3]&=e_1-e_2,\\
[O_4]&=e_1+e_2-e_3,&[O_5]&=e_1-e_2,&[O_6]&=e_1+e_2+e_3,\\
[O_7]&=2e_1-e_2-e_3,&[O_8]&=2e_1+2e_2-e_3,&[O_9]&=e_1-2e_2+e_3,
\end{align*}
that is, by
\[e_2,\quad e_3,\quad e_1+e_2-e_3,\quad
2e_1-e_2-e_3,\quad e_1-2e_2+e_3\]
(since $3[O_3]=[O_7]+[O_9]$).
Thus $\alpha(V)$ is the volume of the domain given
by
\[\begin{cases}
x=1,\,y>0,\,z>0,\\
x+y-z>0,\\
2x-y-z>0,\\
x-2y+z>0.
\end{cases}\]
Using the description above, $\alpha(V)$
is the volume of
\[\begin{cases}
0<y,\,0<z,\\
z-y<1,\\
y+z<2,\\
2y-z<1.
\end{cases}\qquad
\lower 1.3cm\hbox{
\begin{pspicture}(-0.5,-0.5)(2.5,2.5)
\psset{unit=0.75cm}
\rput[bl](2.1,-0.5){$y$}
\rput[bl](-0.5,2.1){$z$}
\psline[linewidth=0.01](1,-0.5)(1,2.5)
\psline[linewidth=0.01](2,-0.5)(2,2.5)
\psline[linewidth=0.01](-0.5,1)(2.5,1)
\psline[linewidth=0.01](-0.5,2)(2.5,2)
\pspolygon[fillstyle=hlines,hatchsep=3pt,%
hatchangle=-45,linestyle=none](0,0)(0.5,0)(1,1)(0.5,1.5)(0,1)(0,0)
\psline{->}(-0.5,0)(2.5,0)
\psline{->}(0,-0.5)(0,2.5)
\psline(1,-0.1)(1,0.1)
\psline(2,-0.1)(2,0.1)
\psline(-0.1,1)(0.1,1)
\psline(-0.1,2)(0.1,2)
\psline(0.25,-0.5)(1.55,2.1)
\psline(-0.5,0.5)(1.1,2.1)
\psline(-0.1,2.1)(2.1,-0.1)
\end{pspicture}}
\]
Therefore $\alpha(V)=1$.
\par
If $V$ is given by the equation \eqref{equ:alpha:rkfour},
then $\Gal(K/\QQ)=\ZZ/2\ZZ$ and the orbits of the Galois
action on the $27$ lines are given by
\begin{align*}
O_1&=\{L(0)\},&O_2&=\{L(1),L(2)\},&O_3&=\{L'(0),L''(0)\},\\
O_4&=\{L'(1),L''(2)\},&O_5&=\{L'(2),L''(1)\},\\
O_6&=\{M(0),M'(0)\},&O_7&=\{M(1),M'(2)\},&O_8&=\{M(2),M'(1)\},\\
O_9&=\{M''(0)\},&O_{10}&=\{M''(1),M''(2)\},\\
O_{11}&=\{N(0),N''(0)\},&O_{12}&=\{N(1),N''(2)\},&O_{13}&=\{N(2),N''(1)\},\\
O_{14}&=\{N'(0)\},&O_{15}&=\{N'(1),N'(2)\}.
\end{align*}
A basis of the Picard group is given by
\[e_1=[\Lambda]-[E_5],\quad e_2=[E_1],\quad e_3=[E_2]+[E_3],
\quad e_4=[E_4]-2[E_5]+[E_6].\]
The effective cone $\Ceff(V)$ is generated by
\begin{align*}
[O_1]&=e_2,&[O_2]&=e_3,\\
[O_3]&=3e_1-2e_2-e_3-e_4,&[O_4]&=3e_1-2e_3-e_4,\\
[O_5]&=3e_1-2e_2-e_3-e_4,&[O_6]&=e_1,\\
[O_7]&=e_1-e_4,&[O_8]&=e_1,\\
[O_9]&=2e_1-e_2-e_3-e_4,&[O_{10}]&=4e_1-2e_2-2e_3-e_4,\\
[O_{11}]&=2e_1-e_3-e_4,&[O_{12}]&=2e_1-e_3-e_4,\\
[O_{13}]&=2e_1-e_3,&[O_{14}]&=e_1-e_2,\\
[O_{15}]&=2e_1-2e_2-e_4.
\end{align*}
Since $[O_3]=[O_5]=[O_9]+[O_{14}]$ and $[O_{11}]=
[O_{12}]=[O_9]+[O_2]$,
we get that $\Ceff(V)$ is generated by
\[\begin{split}
e_2,\quad e_3,\quad 3e_1-2e_3-e_4,\quad e_1-e_4,\quad
2e_1-e_2-e_3-e_4,\\
4e_1-2e_2-2e_3-e_4,\quad 2e_1-e_3,\quad e_1-e_2,\quad
2e_1-2e_2-e_4.
\end{split}\]
The anticanonical class is given by
\[\antican=3e_1-e_2-e_3-e_4.\]
Thus $\alpha(V)$ is the volume of the domain
\[\begin{cases}
3x-y-z-t=1,\\
y>0,\quad z>0,\\
x-y>0,\\
2x-z>0,\\
x-t>0,\\
3x-2z-t>0,\\
2x-y-z-t>0,\\
4x-2y-2z-t>0,\\
2x-2y-t>0,
\end{cases}\]
that is, of the domain $P$ in $\RR^3$ given by
\[\begin{cases}
y>0,\quad z>0,\\
x-y>0,\\
2x-z>0,\\
1-2x+y+z>0,\\
1+y-z>0,\\
1-x>0,\\
1+x-y-z>0,\\
1-x-y+z>0.
\end{cases}\]
We compute its volume as follows:
decompose $P$ into cones with appex $(0,0,0)$ and supported
by the faces not containing this point. Thus we consider
the following faces of $P$:
\begin{align*}
F_1&:\quad1-x=0,&F_2&:\quad1-2x+y+z=0,\\
F_3&:\quad1+y-z=0,&F_4&:\quad 1+x-y-z=0,\\
F_5&:\quad1-x-y+z=0.
\end{align*}
One has
\[\alpha(V)=\Vol(P)=\frac 13\sum_{i=1}^5\Area(F_i).\]
The area of $F_1$ is the volume of the domain
\[\begin{cases}
y>0,\quad z>0,\\
1-y>0,\\
2-z>0,\\
-1+y+z>0,\\
1+y-z>0,\\
2-y-z>0,\\
z-y>0,
\end{cases}\qquad
\lower 1.3cm\hbox{
\begin{pspicture}(-0.5,-0.5)(2.5,2.5)
\psset{unit=0.75cm}
\rput[bl](2.1,-0.5){$y$}
\rput[bl](-0.5,2.1){$z$}
\psline[linewidth=0.01](1,-0.5)(1,2.5)
\psline[linewidth=0.01](2,-0.5)(2,2.5)
\psline[linewidth=0.01](-0.5,1)(2.5,1)
\psline[linewidth=0.01](-0.5,2)(2.5,2)
\pspolygon[fillstyle=hlines,hatchsep=3pt,%
hatchangle=-45,linestyle=none](0,1)(0.5,0.5)(1,1)(0.5,1.5)(0,1)
\psline{->}(-0.5,0)(2.5,0)
\psline{->}(0,-0.5)(0,2.5)
\psline(-0.1,-0.1)(2.1,2.1)
\psline(-0.1,1.1)(1.1,-0.1)
\psline(1,-0.5)(1,2.5)
\psline(-0.1,0.9)(1.1,2.1)
\psline(-0.1,2.1)(2.1,-0.1)
\psline(-0.5,2)(2.5,2)
\end{pspicture}}\]
and we get $\Area(F_1)=\frac 12$. For $F_2$ we have the equations
\[\begin{cases}
y>0,\\
-1+2x-y>0,\\
x-y>0,\\
1+y>0,\\
2-2x+2y>0,\\
1-x>0,\\
2-x>0,\\
x-2y>0.
\end{cases}\qquad
\lower 1.3cm
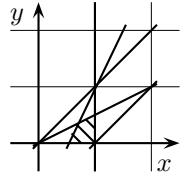
\hbox{
\begin{pspicture}(-0.5,-0.5)(2.5,2.5)
\psset{unit=0.75cm}
\rput[bl](2.1,-0.5){$x$}
\rput[bl](-0.5,2.1){$y$}
\psline[linewidth=0.01](1,-0.5)(1,2.5)
\psline[linewidth=0.01](2,-0.5)(2,2.5)
\psline[linewidth=0.01](-0.5,1)(2.5,1)
\psline[linewidth=0.01](-0.5,2)(2.5,2)
\pspolygon[fillstyle=hlines,hatchsep=3pt,%
hatchangle=-45,linestyle=none](0.5,0)(1,0)(1,0.5)(0.666666,0.333333)(0.5,0)
\psline{->}(-0.5,0)(2.5,0)
\psline{->}(0,-0.5)(0,2.5)
\psline(-0.1,-0.1)(2.1,2.1)
\psline(1,-0.5)(1,2.5)
\psline(0.9,-0.1)(2.1,1.1)
\psline(0.495,-0.1)(1.55,2.1)
\psline(-0.1,-0.05)(2.1,1.05)
\end{pspicture}}
\]
We get $\Area(F_2)=\frac 16$. For $F_3$ we have the same
equations and the same area. For $F_4$ we have the equations
\[\begin{cases}
y>0,\\
1+x-y>0,\\
x-y>0,\\
-1+x+y>0,\\
2-x>0,\\
-x+2y>0,\\
1-x>0,\\
2-2y>0.\\
\end{cases}\qquad
\lower 1.3cm
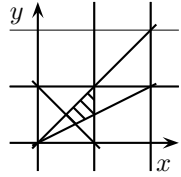
\hbox{
\begin{pspicture}(-0.5,-0.5)(2.5,2.5)
\psset{unit=0.75cm}
\rput[bl](2.1,-0.5){$x$}
\rput[bl](-0.5,2.1){$y$}
\psline[linewidth=0.01](1,-0.5)(1,2.5)
\psline[linewidth=0.01](2,-0.5)(2,2.5)
\psline[linewidth=0.01](-0.5,1)(2.5,1)
\psline[linewidth=0.01](-0.5,2)(2.5,2)
\pspolygon[fillstyle=hlines,hatchsep=3pt,%
hatchangle=-45,linestyle=none](0.5,0.5)(0.666666,0.333333)(1,0.5)
(1,1)(0.5,0.5)
\psline{->}(-0.5,0)(2.5,0)
\psline{->}(0,-0.5)(0,2.5)
\psline(-0.1,-0.1)(2.1,2.1)
\psline(-0.1,-0.05)(2.1,1.05)
\psline(-0.1,1.1)(1.1,-0.1)
\psline(1,-0.5)(1,2.5)
\psline(2,-0.5)(2,2.5)
\psline(-0.5,1)(2.5,1)
\end{pspicture}}\]
We find $\Area(F_4)=1/8+1/24=1/6$. The face $F_5$
is given by the same equations and $\Area(F_5)=1/6$.
Finally
\[\alpha(V)=\frac 13\Bigl(\frac 12+\frac 46\Bigr)=\frac 7{18}.\qed\]
\noqed
\end{proof}


\section{Some statistical formulae}
\label{section:stats}
 
With the program of D. J. Bernstein \cite{bernstein:enumerating}, given 
a number $B$, we can compute the value of $\nUH(2^r)$
for $100\varleq 2^r\varleq B$. Thus we get
a family of pairs 
$(B_i,\nUH(B_i))_{1\varleq i\varleq N}$, where
$B_i$ are consecutive powers of $2$ for $i<N$ and $B_N=B$.
For any $i$ between $1$ and $N$, let
\[x_i=\log(B_i)\quad\text{and}\quad y_i=\nUH(B_i)/B_i.\]
We expect an asymptotic of the form
\[\nUH(B)=BP(\log(B))+o(B),\]
where $P$ is a polynomial of degree $t-1$ with a dominant
coefficient equal to $\thetaH(V)$.
Thus we look for a polynomial $Q$
of degree $t-1$ such that
\[\sum_{i=1}^N(Q(x_i)-y_i)^2\]
is minimal and compare the
leading coefficient of $Q$ with $\thetaH(V)$.

\begin{notas}
Let $R(X,Y)$ be a polynomial in $\QQ[X,Y]$ and denote by
$\langle R(X,Y)\rangle$ the mean value of
$(R(x_i,y_i))_{1\varleq i\varleq N}$, that is,
\[\langle R(X,Y)\rangle=\frac 1N\sum_{i=1}^NR(x_i,y_i).\]
If $t=2$ the leading coefficient of $Q$ (if it is uniquely defined)
is given by
\[A_1=\frac{\langle XY\rangle-\langle Y\rangle\langle X\rangle}{
\langle X^2\rangle-\langle X\rangle^2}.\]
If $t=3$ the leading coefficient is
\[A_2=\frac{\langle YX^2\rangle-\langle Y\rangle\langle X^2\rangle
-\frac{(\langle X^3\rangle-\langle X\rangle\langle X^2\rangle)
(\langle YX\rangle-\langle Y\rangle\langle X\rangle)}{
\langle X^2\rangle-\langle X\rangle^2}}{
\langle X^4\rangle-\langle X^2\rangle^2-
\frac{(\langle X^3\rangle-\langle X\rangle
\langle X^2\rangle)^2}{\langle X^2\rangle-\langle X\rangle^2}}.\]
If $t=4$, the leading coefficient is
\[A_3=\frac{\langle YX^3\rangle-\langle Y\rangle\langle X^3\rangle
-\frac{(\langle X^4\rangle-\langle X\rangle\langle X^3\rangle)
(\langle YX\rangle-\langle Y\rangle\langle X\rangle)}{
\langle X^2\rangle-\langle X\rangle^2}-\frac{\beta\delta}{\gamma}}{
\langle X^6\rangle-\langle X^3\rangle^2-
\frac{(\langle X^4\rangle-\langle X\rangle
\langle X^3\rangle)^2}{\langle X^2\rangle-\langle X\rangle^2}-
\frac{\beta^2}{\gamma}},\]
with
\begin{align*}
\beta&=\langle X^5\rangle-\langle X^3\rangle\langle X^2\rangle
-\frac{\langle X^3\rangle-\langle X\rangle\langle X^2\rangle}{
\langle X^2\rangle-\langle X\rangle^2}(\langle X^4\rangle
-\langle X^3\rangle\langle X\rangle),\\
\gamma&=\langle X^4\rangle-\langle X^2\rangle^2
-\frac{(\langle X^3\rangle-\langle X\rangle\langle X^2\rangle)^2}{
\langle X^2\rangle-\langle X\rangle^2},\\
\delta&=\langle YX^2\rangle-\langle Y\rangle\langle X^2\rangle
-\frac{\langle X^3\rangle-\langle X\rangle\langle X^2\rangle}{
\langle X^2\rangle-\langle X\rangle^2}(\langle YX\rangle
-\langle Y\rangle\langle X\rangle).
\end{align*}
In the next section, we denote by 
$\thetaH^{\text{stat}}(V)$ the leading
coefficient $A_{t-1}$.
\end{notas}

\section{Presentation of the results}
\label{section:results}

We consider only cubic surfaces of the form
\eqref{equ:alpha:rktwo}, \eqref{equ:alpha:rkthree},
or \eqref{equ:alpha:rkfour}.
By \cite[Lemme 1]{colliotkanevskysansuc:cubic}, the corresponding
surface $V$ is $\QQ$-rational and, in particular, $\Br(V)=0$.
Thus the Brauer-Manin obstruction to weak approximation is void and
\[V(\Adeles_\QQ)^{\Br}=V(\Adeles_\QQ)=
\prod_{v\in\Val(\QQ)}V(\QQ_v).\]
Moreover,
\[\beta(V)=\cardinal H^1(\QQ,\Pic(\Vbar))=1.\]
By \eqref{equ:constant:tamagawa} and \eqref{equ:constant:constant}, the
constant $\thetaH(V)$ may be written as
\[\thetaH(V)=\alpha(V)\omegaH(V(\Adeles_\QQ))\]
Using remark \ref{rem:euler:euler} we get 
\begin{multline*}
\thetaH(V)=\alpha(V)\lim_{s\to 1}(s-1)^{t+2}\zeta_E(s)\times
\omegaHof\infty(V(\RR))\\
\times\prod_{p\divise 3abcd}\lambda_p
\omegaHp(V(\QQ_p))\times\prod_{i=0}^3C_i,
\end{multline*}
where $E$ is the \'etale algebra defined in \ref{notas:lines:etale}.
The residue of the zeta function could have been computed directly
(see, for example, \cite[chapter 4]{cohen:computational}), but
instead we used PARI.
The volume at the real place is given by the formula
\[\frac 12\int_{\left\{(x,y,z,t)\left\vert\left\{
\setbox\strutbox\hbox{\vphantom{$\scriptstyle ($}}\!\!\!\!
\begin{array}{l}
\scriptstyle ax^3+by^3+cz^3+dt^3=0\\
\scriptstyle\sup(\vert x\vert,\vert y\vert,\vert z\vert,\vert t\vert)\leq 1
\end{array}\!\!\!\!
\right.\right.\right\}}
\boldsymbol\omega_{L}(x,y,z,t),\]
where $\boldsymbol\omega_L$ is the Leray form 
\[\boldsymbol\omega_L(x,y,z,t)=\frac{\sqrt[3]{d}}{3
{(ax^3+by^3+cz^3)^{2/3}}}\Haar x\Haar y\Haar z.\]
Decomposing the domain of integration (and using the various
expressions of the Leray form) it is possible to remove the singularities
of this integral which is then easily estimated on a computer.
The factors corresponding to the bad places have been described in
section \ref{section:bad} and the constants $C_0$, $C_1$, $C_2$,
and $C_3$ may be computed directly as in section
\ref{section:euler}.
\par
We considered the following examples: for the
cubic surfaces with a Picard group of rank $2$ we used
\begin{align}
\tag{$S_1$}X^3+Y^3+2Z^3+4T^3&=0,\\
\tag{$S_2$}X^3+Y^3+5Z^3+25T^3&=0,\\
\tag{$S_3$}X^3+Y^3+3Z^3+9T^3&=0.\\
\noalign{\noindent For the rank $3$ case:}
\tag{$S_4$}X^3+Y^3+2Z^3+2T^3&=0,\\
\tag{$S_5$}X^3+Y^3+5Z^3+5T^3&=0,\\
\tag{$S_6$}X^3+Y^3+7Z^3+7T^3&=0,\\
\tag{$S_7$}2X^3+2Y^3+3Z^3+3T^3&=0,\\
\noalign{\noindent and for rank $4$:}
\tag{$S_8$}X^3+Y^3+Z^3+T^3&=0.
\end{align}
We draw below the corresponding experimental curves in which we compare
the value of $\nUH(B)/(B(\log B)^{t-1})$ with $\thetaH(V)$.
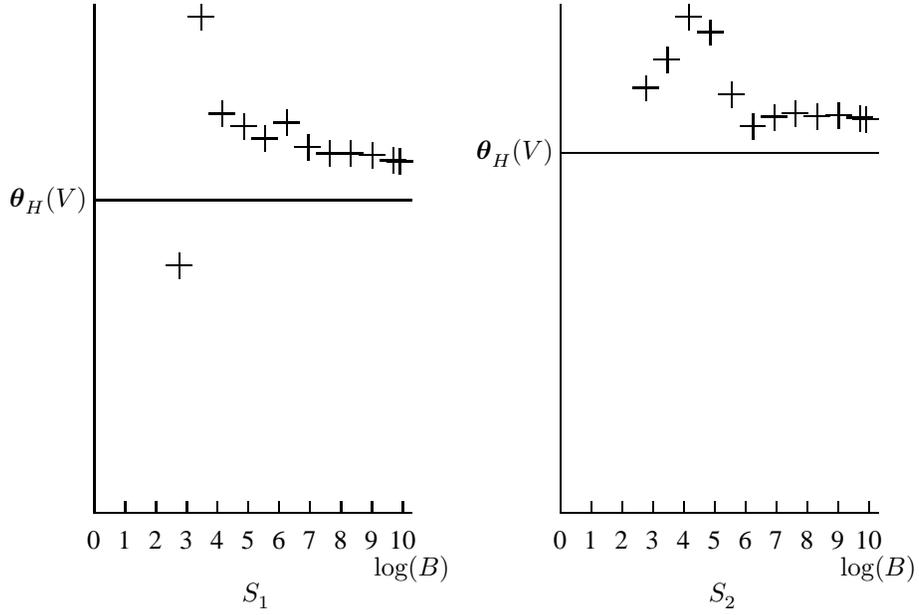
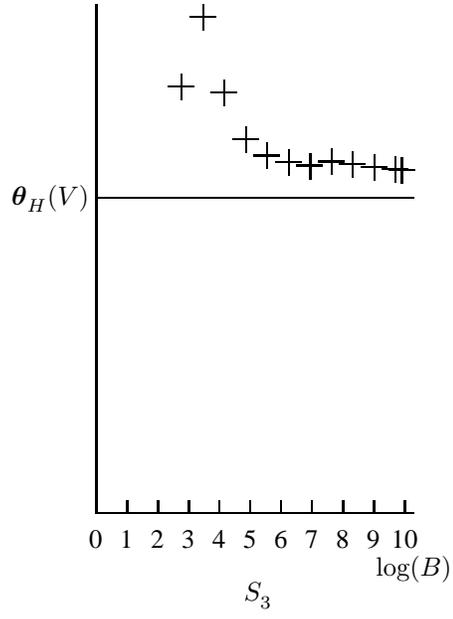
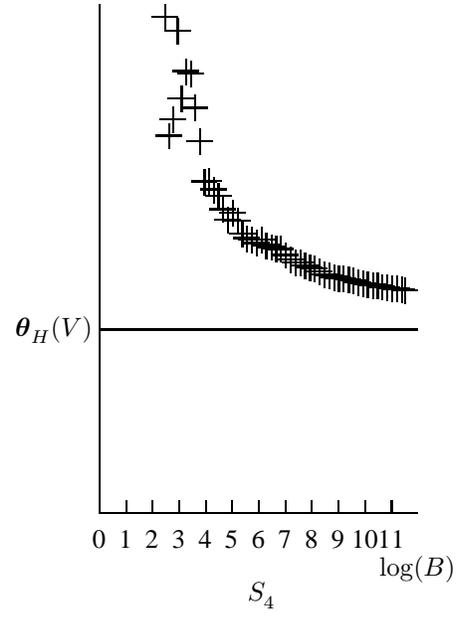
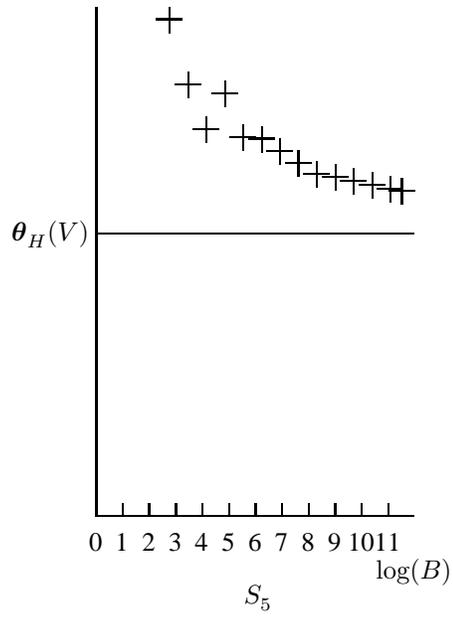
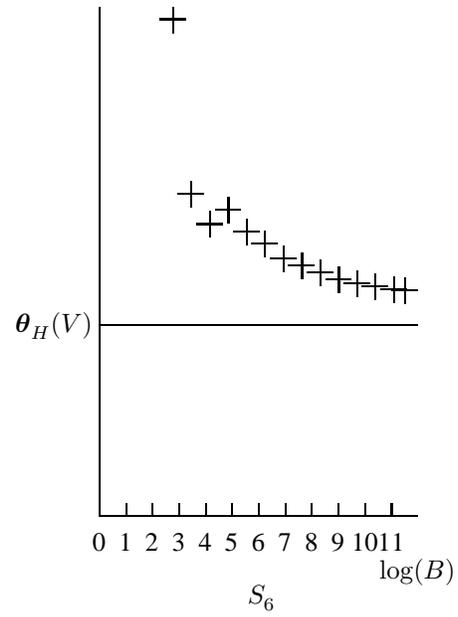
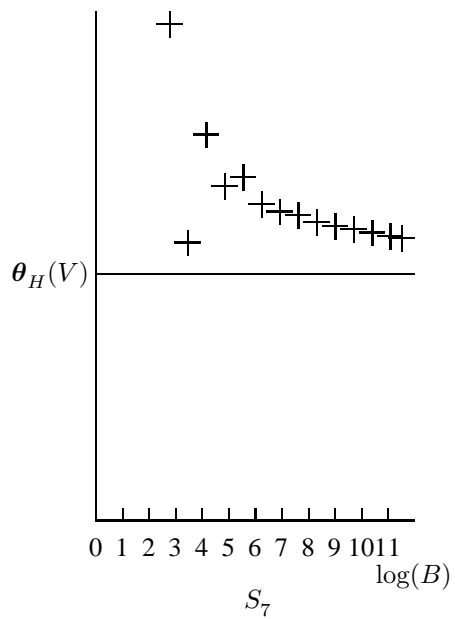
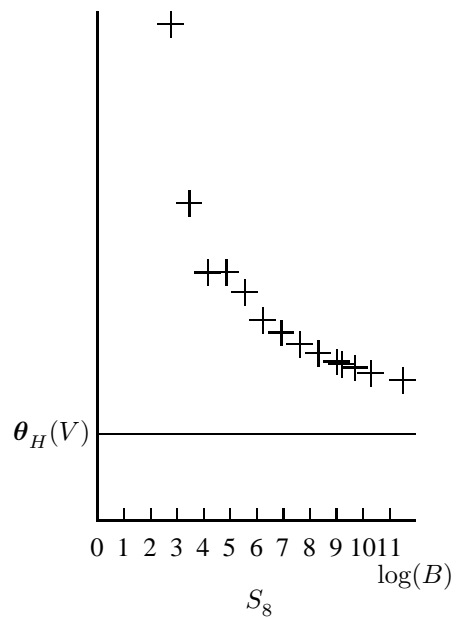
\[\vbox{%
\hsize\graphsize
\hbox to\hsize{\hfil
{\setlength{\unitlength}{0.240900pt}
\begin{picture}(600,900)(0,0)
\put(50.000,88.000){\rule[-0.200pt]{120.450pt}{0.400pt}}\put(50.000,88.000){\rule[-0.200pt]{0.400pt}{192.720pt}}\put(50.000,88.000){\rule[-0.200pt]{0.400pt}{4.818pt}}
\put(50.000,47.0){\makebox(0,0)[c]{0}}
\put(98.468,88.000){\rule[-0.200pt]{0.400pt}{4.818pt}}
\put(98.468,47.0){\makebox(0,0)[c]{1}}
\put(146.936,88.000){\rule[-0.200pt]{0.400pt}{4.818pt}}
\put(146.936,47.0){\makebox(0,0)[c]{2}}
\put(195.404,88.000){\rule[-0.200pt]{0.400pt}{4.818pt}}
\put(195.404,47.0){\makebox(0,0)[c]{3}}
\put(243.872,88.000){\rule[-0.200pt]{0.400pt}{4.818pt}}
\put(243.872,47.0){\makebox(0,0)[c]{4}}
\put(292.340,88.000){\rule[-0.200pt]{0.400pt}{4.818pt}}
\put(292.340,47.0){\makebox(0,0)[c]{5}}
\put(340.808,88.000){\rule[-0.200pt]{0.400pt}{4.818pt}}
\put(340.808,47.0){\makebox(0,0)[c]{6}}
\put(389.276,88.000){\rule[-0.200pt]{0.400pt}{4.818pt}}
\put(389.276,47.0){\makebox(0,0)[c]{7}}
\put(437.744,88.000){\rule[-0.200pt]{0.400pt}{4.818pt}}
\put(437.744,47.0){\makebox(0,0)[c]{8}}
\put(486.212,88.000){\rule[-0.200pt]{0.400pt}{4.818pt}}
\put(486.212,47.0){\makebox(0,0)[c]{9}}
\put(534.680,88.000){\rule[-0.200pt]{0.400pt}{4.818pt}}
\put(534.680,47.0){\makebox(0,0)[c]{10}}
\put(550.000,0.0){\makebox(0,0)[c]{$\log(B)$}}
\put(164.382,478.000){\rule[-0.200pt]{9.636pt}{0.400pt}}
\put(184.382,458.000){\rule[-0.200pt]{0.400pt}{9.636pt}}
\put(197.977,868.000){\rule[-0.200pt]{9.636pt}{0.400pt}}
\put(217.977,848.000){\rule[-0.200pt]{0.400pt}{9.636pt}}
\put(231.573,716.333){\rule[-0.200pt]{9.636pt}{0.400pt}}
\put(251.573,696.333){\rule[-0.200pt]{0.400pt}{9.636pt}}
\put(265.168,696.214){\rule[-0.200pt]{9.636pt}{0.400pt}}
\put(285.168,676.214){\rule[-0.200pt]{0.400pt}{9.636pt}}
\put(298.764,677.062){\rule[-0.200pt]{9.636pt}{0.400pt}}
\put(318.764,657.062){\rule[-0.200pt]{0.400pt}{9.636pt}}
\put(332.359,701.889){\rule[-0.200pt]{9.636pt}{0.400pt}}
\put(352.359,681.889){\rule[-0.200pt]{0.400pt}{9.636pt}}
\put(365.955,663.656){\rule[-0.200pt]{9.636pt}{0.400pt}}
\put(385.955,643.656){\rule[-0.200pt]{0.400pt}{9.636pt}}
\put(399.550,653.611){\rule[-0.200pt]{9.636pt}{0.400pt}}
\put(419.550,633.611){\rule[-0.200pt]{0.400pt}{9.636pt}}
\put(433.146,653.365){\rule[-0.200pt]{9.636pt}{0.400pt}}
\put(453.146,633.365){\rule[-0.200pt]{0.400pt}{9.636pt}}
\put(466.741,650.969){\rule[-0.200pt]{9.636pt}{0.400pt}}
\put(486.741,630.969){\rule[-0.200pt]{0.400pt}{9.636pt}}
\put(500.337,642.422){\rule[-0.200pt]{9.636pt}{0.400pt}}
\put(520.337,622.422){\rule[-0.200pt]{0.400pt}{9.636pt}}
\put(510.000,641.117){\rule[-0.200pt]{9.636pt}{0.400pt}}
\put(530.000,621.117){\rule[-0.200pt]{0.400pt}{9.636pt}}
\put(50.0,580.140){\rule[-0.200pt]{120.450pt}{0.400pt}}
\put(40.0,580.140){\makebox(0,0)[r]{$\thetaH(V)$}}
\end{picture}}
\hfil}
\hbox to\hsize{
\hfil$S_1$\hfil}}
\vbox{%
\hsize\graphsize
\hbox to\hsize{\hfil
{\setlength{\unitlength}{0.240900pt}
\begin{picture}(600,900)(0,0)
\put(50.000,88.000){\rule[-0.200pt]{120.450pt}{0.400pt}}\put(50.000,88.000){\rule[-0.200pt]{0.400pt}{192.720pt}}\put(50.000,88.000){\rule[-0.200pt]{0.400pt}{4.818pt}}
\put(50.000,47.0){\makebox(0,0)[c]{0}}
\put(98.468,88.000){\rule[-0.200pt]{0.400pt}{4.818pt}}
\put(98.468,47.0){\makebox(0,0)[c]{1}}
\put(146.936,88.000){\rule[-0.200pt]{0.400pt}{4.818pt}}
\put(146.936,47.0){\makebox(0,0)[c]{2}}
\put(195.404,88.000){\rule[-0.200pt]{0.400pt}{4.818pt}}
\put(195.404,47.0){\makebox(0,0)[c]{3}}
\put(243.872,88.000){\rule[-0.200pt]{0.400pt}{4.818pt}}
\put(243.872,47.0){\makebox(0,0)[c]{4}}
\put(292.340,88.000){\rule[-0.200pt]{0.400pt}{4.818pt}}
\put(292.340,47.0){\makebox(0,0)[c]{5}}
\put(340.808,88.000){\rule[-0.200pt]{0.400pt}{4.818pt}}
\put(340.808,47.0){\makebox(0,0)[c]{6}}
\put(389.276,88.000){\rule[-0.200pt]{0.400pt}{4.818pt}}
\put(389.276,47.0){\makebox(0,0)[c]{7}}
\put(437.744,88.000){\rule[-0.200pt]{0.400pt}{4.818pt}}
\put(437.744,47.0){\makebox(0,0)[c]{8}}
\put(486.212,88.000){\rule[-0.200pt]{0.400pt}{4.818pt}}
\put(486.212,47.0){\makebox(0,0)[c]{9}}
\put(534.680,88.000){\rule[-0.200pt]{0.400pt}{4.818pt}}
\put(534.680,47.0){\makebox(0,0)[c]{10}}
\put(550.000,0.0){\makebox(0,0)[c]{$\log(B)$}}
\put(164.382,756.571){\rule[-0.200pt]{9.636pt}{0.400pt}}
\put(184.382,736.571){\rule[-0.200pt]{0.400pt}{9.636pt}}
\put(197.977,801.143){\rule[-0.200pt]{9.636pt}{0.400pt}}
\put(217.977,781.143){\rule[-0.200pt]{0.400pt}{9.636pt}}
\put(231.573,868.000){\rule[-0.200pt]{9.636pt}{0.400pt}}
\put(251.573,848.000){\rule[-0.200pt]{0.400pt}{9.636pt}}
\put(265.168,844.122){\rule[-0.200pt]{9.636pt}{0.400pt}}
\put(285.168,824.122){\rule[-0.200pt]{0.400pt}{9.636pt}}
\put(298.764,746.125){\rule[-0.200pt]{9.636pt}{0.400pt}}
\put(318.764,726.125){\rule[-0.200pt]{0.400pt}{9.636pt}}
\put(332.359,696.214){\rule[-0.200pt]{9.636pt}{0.400pt}}
\put(352.359,676.214){\rule[-0.200pt]{0.400pt}{9.636pt}}
\put(365.955,711.304){\rule[-0.200pt]{9.636pt}{0.400pt}}
\put(385.955,691.304){\rule[-0.200pt]{0.400pt}{9.636pt}}
\put(399.550,716.685){\rule[-0.200pt]{9.636pt}{0.400pt}}
\put(419.550,696.685){\rule[-0.200pt]{0.400pt}{9.636pt}}
\put(433.146,711.884){\rule[-0.200pt]{9.636pt}{0.400pt}}
\put(453.146,691.884){\rule[-0.200pt]{0.400pt}{9.636pt}}
\put(466.741,713.513){\rule[-0.200pt]{9.636pt}{0.400pt}}
\put(486.741,693.513){\rule[-0.200pt]{0.400pt}{9.636pt}}
\put(500.337,709.376){\rule[-0.200pt]{9.636pt}{0.400pt}}
\put(520.337,689.376){\rule[-0.200pt]{0.400pt}{9.636pt}}
\put(510.000,707.055){\rule[-0.200pt]{9.636pt}{0.400pt}}
\put(530.000,687.055){\rule[-0.200pt]{0.400pt}{9.636pt}}
\put(50.0,654.179){\rule[-0.200pt]{120.450pt}{0.400pt}}
\put(40.0,654.179){\makebox(0,0)[r]{$\thetaH(V)$}}
\end{picture}}

\hfil}
\hbox to\hsize{
\hfil$S_2$\hfil}}\]
\vfil
\penalty 600
\vfilneg
\[\vbox{%
\hsize\graphsize
\hbox to\hsize{\hfil
{\setlength{\unitlength}{0.240900pt}
\begin{picture}(600,900)(0,0)
\put(50.000,88.000){\rule[-0.200pt]{120.450pt}{0.400pt}}\put(50.000,88.000){\rule[-0.200pt]{0.400pt}{192.720pt}}\put(50.000,88.000){\rule[-0.200pt]{0.400pt}{4.818pt}}
\put(50.000,47.0){\makebox(0,0)[c]{0}}
\put(98.468,88.000){\rule[-0.200pt]{0.400pt}{4.818pt}}
\put(98.468,47.0){\makebox(0,0)[c]{1}}
\put(146.936,88.000){\rule[-0.200pt]{0.400pt}{4.818pt}}
\put(146.936,47.0){\makebox(0,0)[c]{2}}
\put(195.403,88.000){\rule[-0.200pt]{0.400pt}{4.818pt}}
\put(195.403,47.0){\makebox(0,0)[c]{3}}
\put(243.871,88.000){\rule[-0.200pt]{0.400pt}{4.818pt}}
\put(243.871,47.0){\makebox(0,0)[c]{4}}
\put(292.339,88.000){\rule[-0.200pt]{0.400pt}{4.818pt}}
\put(292.339,47.0){\makebox(0,0)[c]{5}}
\put(340.807,88.000){\rule[-0.200pt]{0.400pt}{4.818pt}}
\put(340.807,47.0){\makebox(0,0)[c]{6}}
\put(389.274,88.000){\rule[-0.200pt]{0.400pt}{4.818pt}}
\put(389.274,47.0){\makebox(0,0)[c]{7}}
\put(437.742,88.000){\rule[-0.200pt]{0.400pt}{4.818pt}}
\put(437.742,47.0){\makebox(0,0)[c]{8}}
\put(486.210,88.000){\rule[-0.200pt]{0.400pt}{4.818pt}}
\put(486.210,47.0){\makebox(0,0)[c]{9}}
\put(534.678,88.000){\rule[-0.200pt]{0.400pt}{4.818pt}}
\put(534.678,47.0){\makebox(0,0)[c]{10}}
\put(550.000,0.0){\makebox(0,0)[c]{$\log(B)$}}
\put(164.381,758.313){\rule[-0.200pt]{9.636pt}{0.400pt}}
\put(184.381,738.313){\rule[-0.200pt]{0.400pt}{9.636pt}}
\put(197.977,868.000){\rule[-0.200pt]{9.636pt}{0.400pt}}
\put(217.977,848.000){\rule[-0.200pt]{0.400pt}{9.636pt}}
\put(231.572,748.156){\rule[-0.200pt]{9.636pt}{0.400pt}}
\put(251.572,728.156){\rule[-0.200pt]{0.400pt}{9.636pt}}
\put(265.167,675.612){\rule[-0.200pt]{9.636pt}{0.400pt}}
\put(285.167,655.612){\rule[-0.200pt]{0.400pt}{9.636pt}}
\put(298.762,649.768){\rule[-0.200pt]{9.636pt}{0.400pt}}
\put(318.762,629.768){\rule[-0.200pt]{0.400pt}{9.636pt}}
\put(332.358,639.823){\rule[-0.200pt]{9.636pt}{0.400pt}}
\put(352.358,619.823){\rule[-0.200pt]{0.400pt}{9.636pt}}
\put(365.953,633.771){\rule[-0.200pt]{9.636pt}{0.400pt}}
\put(385.953,613.771){\rule[-0.200pt]{0.400pt}{9.636pt}}
\put(399.548,640.246){\rule[-0.200pt]{9.636pt}{0.400pt}}
\put(419.548,620.246){\rule[-0.200pt]{0.400pt}{9.636pt}}
\put(433.144,636.358){\rule[-0.200pt]{9.636pt}{0.400pt}}
\put(453.144,616.358){\rule[-0.200pt]{0.400pt}{9.636pt}}
\put(466.739,632.080){\rule[-0.200pt]{9.636pt}{0.400pt}}
\put(486.739,612.080){\rule[-0.200pt]{0.400pt}{9.636pt}}
\put(500.334,628.140){\rule[-0.200pt]{9.636pt}{0.400pt}}
\put(520.334,608.140){\rule[-0.200pt]{0.400pt}{9.636pt}}
\put(510.000,626.963){\rule[-0.200pt]{9.636pt}{0.400pt}}
\put(530.000,606.963){\rule[-0.200pt]{0.400pt}{9.636pt}}
\put(50.0,582.819){\rule[-0.200pt]{120.450pt}{0.400pt}}
\put(40.0,582.819){\makebox(0,0)[r]{$\thetaH(V)$}}
\end{picture}}
\hfil}
\hbox to\hsize{
\hfil$S_3$\hfil}}
\vbox{%
\hsize\graphsize
\hbox to\hsize{\hfil
{\setlength{\unitlength}{0.240900pt}
\begin{picture}(600,900)(0,0)
\put(50.000,88.000){\rule[-0.200pt]{120.450pt}{0.400pt}}\put(50.000,88.000){\rule[-0.200pt]{0.400pt}{192.720pt}}\put(50.000,88.000){\rule[-0.200pt]{0.400pt}{4.818pt}}
\put(50.000,47.0){\makebox(0,0)[c]{0}}
\put(91.692,88.000){\rule[-0.200pt]{0.400pt}{4.818pt}}
\put(91.692,47.0){\makebox(0,0)[c]{1}}
\put(133.385,88.000){\rule[-0.200pt]{0.400pt}{4.818pt}}
\put(133.385,47.0){\makebox(0,0)[c]{2}}
\put(175.077,88.000){\rule[-0.200pt]{0.400pt}{4.818pt}}
\put(175.077,47.0){\makebox(0,0)[c]{3}}
\put(216.769,88.000){\rule[-0.200pt]{0.400pt}{4.818pt}}
\put(216.769,47.0){\makebox(0,0)[c]{4}}
\put(258.461,88.000){\rule[-0.200pt]{0.400pt}{4.818pt}}
\put(258.461,47.0){\makebox(0,0)[c]{5}}
\put(300.154,88.000){\rule[-0.200pt]{0.400pt}{4.818pt}}
\put(300.154,47.0){\makebox(0,0)[c]{6}}
\put(341.846,88.000){\rule[-0.200pt]{0.400pt}{4.818pt}}
\put(341.846,47.0){\makebox(0,0)[c]{7}}
\put(383.538,88.000){\rule[-0.200pt]{0.400pt}{4.818pt}}
\put(383.538,47.0){\makebox(0,0)[c]{8}}
\put(425.230,88.000){\rule[-0.200pt]{0.400pt}{4.818pt}}
\put(425.230,47.0){\makebox(0,0)[c]{9}}
\put(466.923,88.000){\rule[-0.200pt]{0.400pt}{4.818pt}}
\put(466.923,47.0){\makebox(0,0)[c]{10}}
\put(508.615,88.000){\rule[-0.200pt]{0.400pt}{4.818pt}}
\put(508.615,47.0){\makebox(0,0)[c]{11}}
\put(550.000,0.0){\makebox(0,0)[c]{$\log(B)$}}
\put(133.601,868.000){\rule[-0.200pt]{9.636pt}{0.400pt}}
\put(153.601,848.000){\rule[-0.200pt]{0.400pt}{9.636pt}}
\put(140.028,680.748){\rule[-0.200pt]{9.636pt}{0.400pt}}
\put(160.028,660.748){\rule[-0.200pt]{0.400pt}{9.636pt}}
\put(145.596,706.289){\rule[-0.200pt]{9.636pt}{0.400pt}}
\put(165.596,686.289){\rule[-0.200pt]{0.400pt}{9.636pt}}
\put(152.760,845.125){\rule[-0.200pt]{9.636pt}{0.400pt}}
\put(172.760,825.125){\rule[-0.200pt]{0.400pt}{9.636pt}}
\put(158.873,739.212){\rule[-0.200pt]{9.636pt}{0.400pt}}
\put(178.873,719.212){\rule[-0.200pt]{0.400pt}{9.636pt}}
\put(165.837,782.355){\rule[-0.200pt]{9.636pt}{0.400pt}}
\put(185.837,762.355){\rule[-0.200pt]{0.400pt}{9.636pt}}
\put(173.171,778.655){\rule[-0.200pt]{9.636pt}{0.400pt}}
\put(193.171,758.655){\rule[-0.200pt]{0.400pt}{9.636pt}}
\put(180.547,724.835){\rule[-0.200pt]{9.636pt}{0.400pt}}
\put(200.547,704.835){\rule[-0.200pt]{0.400pt}{9.636pt}}
\put(187.771,672.157){\rule[-0.200pt]{9.636pt}{0.400pt}}
\put(207.771,652.157){\rule[-0.200pt]{0.400pt}{9.636pt}}
\put(194.736,608.818){\rule[-0.200pt]{9.636pt}{0.400pt}}
\put(214.736,588.818){\rule[-0.200pt]{0.400pt}{9.636pt}}
\put(202.070,609.353){\rule[-0.200pt]{9.636pt}{0.400pt}}
\put(222.070,589.353){\rule[-0.200pt]{0.400pt}{9.636pt}}
\put(209.446,596.146){\rule[-0.200pt]{9.636pt}{0.400pt}}
\put(229.446,576.146){\rule[-0.200pt]{0.400pt}{9.636pt}}
\put(216.670,586.331){\rule[-0.200pt]{9.636pt}{0.400pt}}
\put(236.670,566.331){\rule[-0.200pt]{0.400pt}{9.636pt}}
\put(224.034,565.503){\rule[-0.200pt]{9.636pt}{0.400pt}}
\put(244.034,545.503){\rule[-0.200pt]{0.400pt}{9.636pt}}
\put(231.636,549.377){\rule[-0.200pt]{9.636pt}{0.400pt}}
\put(251.636,529.377){\rule[-0.200pt]{0.400pt}{9.636pt}}
\put(239.182,559.351){\rule[-0.200pt]{9.636pt}{0.400pt}}
\put(259.182,539.351){\rule[-0.200pt]{0.400pt}{9.636pt}}
\put(246.737,547.571){\rule[-0.200pt]{9.636pt}{0.400pt}}
\put(266.737,527.571){\rule[-0.200pt]{0.400pt}{9.636pt}}
\put(254.300,527.281){\rule[-0.200pt]{9.636pt}{0.400pt}}
\put(274.300,507.281){\rule[-0.200pt]{0.400pt}{9.636pt}}
\put(261.837,519.716){\rule[-0.200pt]{9.636pt}{0.400pt}}
\put(281.837,499.716){\rule[-0.200pt]{0.400pt}{9.636pt}}
\put(269.439,517.512){\rule[-0.200pt]{9.636pt}{0.400pt}}
\put(289.439,497.512){\rule[-0.200pt]{0.400pt}{9.636pt}}
\put(276.996,512.318){\rule[-0.200pt]{9.636pt}{0.400pt}}
\put(296.996,492.318){\rule[-0.200pt]{0.400pt}{9.636pt}}
\put(284.523,517.784){\rule[-0.200pt]{9.636pt}{0.400pt}}
\put(304.523,497.784){\rule[-0.200pt]{0.400pt}{9.636pt}}
\put(292.078,508.753){\rule[-0.200pt]{9.636pt}{0.400pt}}
\put(312.078,488.753){\rule[-0.200pt]{0.400pt}{9.636pt}}
\put(299.653,505.904){\rule[-0.200pt]{9.636pt}{0.400pt}}
\put(319.653,485.904){\rule[-0.200pt]{0.400pt}{9.636pt}}
\put(307.211,503.604){\rule[-0.200pt]{9.636pt}{0.400pt}}
\put(327.211,483.604){\rule[-0.200pt]{0.400pt}{9.636pt}}
\put(314.795,502.591){\rule[-0.200pt]{9.636pt}{0.400pt}}
\put(334.795,482.591){\rule[-0.200pt]{0.400pt}{9.636pt}}
\put(322.389,493.108){\rule[-0.200pt]{9.636pt}{0.400pt}}
\put(342.389,473.108){\rule[-0.200pt]{0.400pt}{9.636pt}}
\put(329.984,486.515){\rule[-0.200pt]{9.636pt}{0.400pt}}
\put(349.984,466.515){\rule[-0.200pt]{0.400pt}{9.636pt}}
\put(337.569,481.777){\rule[-0.200pt]{9.636pt}{0.400pt}}
\put(357.569,461.777){\rule[-0.200pt]{0.400pt}{9.636pt}}
\put(345.153,481.556){\rule[-0.200pt]{9.636pt}{0.400pt}}
\put(365.153,461.556){\rule[-0.200pt]{0.400pt}{9.636pt}}
\put(352.744,477.262){\rule[-0.200pt]{9.636pt}{0.400pt}}
\put(372.744,457.262){\rule[-0.200pt]{0.400pt}{9.636pt}}
\put(360.342,474.276){\rule[-0.200pt]{9.636pt}{0.400pt}}
\put(380.342,454.276){\rule[-0.200pt]{0.400pt}{9.636pt}}
\put(367.941,472.408){\rule[-0.200pt]{9.636pt}{0.400pt}}
\put(387.941,452.408){\rule[-0.200pt]{0.400pt}{9.636pt}}
\put(375.536,467.767){\rule[-0.200pt]{9.636pt}{0.400pt}}
\put(395.536,447.767){\rule[-0.200pt]{0.400pt}{9.636pt}}
\put(383.138,463.191){\rule[-0.200pt]{9.636pt}{0.400pt}}
\put(403.138,443.191){\rule[-0.200pt]{0.400pt}{9.636pt}}
\put(390.739,461.304){\rule[-0.200pt]{9.636pt}{0.400pt}}
\put(410.739,441.304){\rule[-0.200pt]{0.400pt}{9.636pt}}
\put(398.336,458.796){\rule[-0.200pt]{9.636pt}{0.400pt}}
\put(418.336,438.796){\rule[-0.200pt]{0.400pt}{9.636pt}}
\put(405.934,457.111){\rule[-0.200pt]{9.636pt}{0.400pt}}
\put(425.934,437.111){\rule[-0.200pt]{0.400pt}{9.636pt}}
\put(413.535,456.226){\rule[-0.200pt]{9.636pt}{0.400pt}}
\put(433.535,436.226){\rule[-0.200pt]{0.400pt}{9.636pt}}
\put(421.133,454.995){\rule[-0.200pt]{9.636pt}{0.400pt}}
\put(441.133,434.995){\rule[-0.200pt]{0.400pt}{9.636pt}}
\put(428.734,453.223){\rule[-0.200pt]{9.636pt}{0.400pt}}
\put(448.734,433.223){\rule[-0.200pt]{0.400pt}{9.636pt}}
\put(436.333,450.563){\rule[-0.200pt]{9.636pt}{0.400pt}}
\put(456.333,430.563){\rule[-0.200pt]{0.400pt}{9.636pt}}
\put(443.934,448.795){\rule[-0.200pt]{9.636pt}{0.400pt}}
\put(463.934,428.795){\rule[-0.200pt]{0.400pt}{9.636pt}}
\put(451.535,447.637){\rule[-0.200pt]{9.636pt}{0.400pt}}
\put(471.535,427.637){\rule[-0.200pt]{0.400pt}{9.636pt}}
\put(459.135,446.529){\rule[-0.200pt]{9.636pt}{0.400pt}}
\put(479.135,426.529){\rule[-0.200pt]{0.400pt}{9.636pt}}
\put(466.736,445.165){\rule[-0.200pt]{9.636pt}{0.400pt}}
\put(486.736,425.165){\rule[-0.200pt]{0.400pt}{9.636pt}}
\put(474.337,443.587){\rule[-0.200pt]{9.636pt}{0.400pt}}
\put(494.337,423.587){\rule[-0.200pt]{0.400pt}{9.636pt}}
\put(481.938,442.610){\rule[-0.200pt]{9.636pt}{0.400pt}}
\put(501.938,422.610){\rule[-0.200pt]{0.400pt}{9.636pt}}
\put(489.539,441.538){\rule[-0.200pt]{9.636pt}{0.400pt}}
\put(509.539,421.538){\rule[-0.200pt]{0.400pt}{9.636pt}}
\put(497.140,440.174){\rule[-0.200pt]{9.636pt}{0.400pt}}
\put(517.140,420.174){\rule[-0.200pt]{0.400pt}{9.636pt}}
\put(504.741,439.205){\rule[-0.200pt]{9.636pt}{0.400pt}}
\put(524.741,419.205){\rule[-0.200pt]{0.400pt}{9.636pt}}
\put(510.000,438.116){\rule[-0.200pt]{9.636pt}{0.400pt}}
\put(530.000,418.116){\rule[-0.200pt]{0.400pt}{9.636pt}}
\put(50.0,376.339){\rule[-0.200pt]{120.450pt}{0.400pt}}
\put(40.0,376.339){\makebox(0,0)[r]{$\thetaH(V)$}}
\end{picture}}
\hfil}
\hbox to\hsize{
\hfil$S_4$\hfil}}\]
\vskip 1cm plus 1cm
\[\vbox{%
\hsize\graphsize
\hbox to\hsize{\hfil
{\setlength{\unitlength}{0.240900pt}
\begin{picture}(600,900)(0,0)
\put(50.000,88.000){\rule[-0.200pt]{120.450pt}{0.400pt}}\put(50.000,88.000){\rule[-0.200pt]{0.400pt}{192.720pt}}\put(50.000,88.000){\rule[-0.200pt]{0.400pt}{4.818pt}}
\put(50.000,47.0){\makebox(0,0)[c]{0}}
\put(91.692,88.000){\rule[-0.200pt]{0.400pt}{4.818pt}}
\put(91.692,47.0){\makebox(0,0)[c]{1}}
\put(133.385,88.000){\rule[-0.200pt]{0.400pt}{4.818pt}}
\put(133.385,47.0){\makebox(0,0)[c]{2}}
\put(175.077,88.000){\rule[-0.200pt]{0.400pt}{4.818pt}}
\put(175.077,47.0){\makebox(0,0)[c]{3}}
\put(216.769,88.000){\rule[-0.200pt]{0.400pt}{4.818pt}}
\put(216.769,47.0){\makebox(0,0)[c]{4}}
\put(258.461,88.000){\rule[-0.200pt]{0.400pt}{4.818pt}}
\put(258.461,47.0){\makebox(0,0)[c]{5}}
\put(300.154,88.000){\rule[-0.200pt]{0.400pt}{4.818pt}}
\put(300.154,47.0){\makebox(0,0)[c]{6}}
\put(341.846,88.000){\rule[-0.200pt]{0.400pt}{4.818pt}}
\put(341.846,47.0){\makebox(0,0)[c]{7}}
\put(383.538,88.000){\rule[-0.200pt]{0.400pt}{4.818pt}}
\put(383.538,47.0){\makebox(0,0)[c]{8}}
\put(425.230,88.000){\rule[-0.200pt]{0.400pt}{4.818pt}}
\put(425.230,47.0){\makebox(0,0)[c]{9}}
\put(466.923,88.000){\rule[-0.200pt]{0.400pt}{4.818pt}}
\put(466.923,47.0){\makebox(0,0)[c]{10}}
\put(508.615,88.000){\rule[-0.200pt]{0.400pt}{4.818pt}}
\put(508.615,47.0){\makebox(0,0)[c]{11}}
\put(550.000,0.0){\makebox(0,0)[c]{$\log(B)$}}
\put(145.596,868.000){\rule[-0.200pt]{9.636pt}{0.400pt}}
\put(165.596,848.000){\rule[-0.200pt]{0.400pt}{9.636pt}}
\put(174.494,765.486){\rule[-0.200pt]{9.636pt}{0.400pt}}
\put(194.494,745.486){\rule[-0.200pt]{0.400pt}{9.636pt}}
\put(203.393,694.667){\rule[-0.200pt]{9.636pt}{0.400pt}}
\put(223.393,674.667){\rule[-0.200pt]{0.400pt}{9.636pt}}
\put(232.292,752.023){\rule[-0.200pt]{9.636pt}{0.400pt}}
\put(252.292,732.023){\rule[-0.200pt]{0.400pt}{9.636pt}}
\put(261.191,683.446){\rule[-0.200pt]{9.636pt}{0.400pt}}
\put(281.191,663.446){\rule[-0.200pt]{0.400pt}{9.636pt}}
\put(290.090,680.222){\rule[-0.200pt]{9.636pt}{0.400pt}}
\put(310.090,660.222){\rule[-0.200pt]{0.400pt}{9.636pt}}
\put(318.989,661.021){\rule[-0.200pt]{9.636pt}{0.400pt}}
\put(338.989,641.021){\rule[-0.200pt]{0.400pt}{9.636pt}}
\put(347.888,641.574){\rule[-0.200pt]{9.636pt}{0.400pt}}
\put(367.888,621.574){\rule[-0.200pt]{0.400pt}{9.636pt}}
\put(376.787,624.540){\rule[-0.200pt]{9.636pt}{0.400pt}}
\put(396.787,604.540){\rule[-0.200pt]{0.400pt}{9.636pt}}
\put(405.685,620.747){\rule[-0.200pt]{9.636pt}{0.400pt}}
\put(425.685,600.747){\rule[-0.200pt]{0.400pt}{9.636pt}}
\put(434.584,614.043){\rule[-0.200pt]{9.636pt}{0.400pt}}
\put(454.584,594.043){\rule[-0.200pt]{0.400pt}{9.636pt}}
\put(463.483,608.114){\rule[-0.200pt]{9.636pt}{0.400pt}}
\put(483.483,588.114){\rule[-0.200pt]{0.400pt}{9.636pt}}
\put(492.382,601.716){\rule[-0.200pt]{9.636pt}{0.400pt}}
\put(512.382,581.716){\rule[-0.200pt]{0.400pt}{9.636pt}}
\put(510.000,598.915){\rule[-0.200pt]{9.636pt}{0.400pt}}
\put(530.000,578.915){\rule[-0.200pt]{0.400pt}{9.636pt}}
\put(50.0,530.661){\rule[-0.200pt]{120.450pt}{0.400pt}}
\put(40.0,530.661){\makebox(0,0)[r]{$\thetaH(V)$}}
\end{picture}}
\hfil}
\hbox to\hsize{
\hfil$S_5$\hfil}}
\vbox{%
\hsize\graphsize
\hbox to\hsize{\hfil
{\setlength{\unitlength}{0.240900pt}
\begin{picture}(600,900)(0,0)
\put(50.000,88.000){\rule[-0.200pt]{120.450pt}{0.400pt}}\put(50.000,88.000){\rule[-0.200pt]{0.400pt}{192.720pt}}\put(50.000,88.000){\rule[-0.200pt]{0.400pt}{4.818pt}}
\put(50.000,47.0){\makebox(0,0)[c]{0}}
\put(91.692,88.000){\rule[-0.200pt]{0.400pt}{4.818pt}}
\put(91.692,47.0){\makebox(0,0)[c]{1}}
\put(133.385,88.000){\rule[-0.200pt]{0.400pt}{4.818pt}}
\put(133.385,47.0){\makebox(0,0)[c]{2}}
\put(175.077,88.000){\rule[-0.200pt]{0.400pt}{4.818pt}}
\put(175.077,47.0){\makebox(0,0)[c]{3}}
\put(216.769,88.000){\rule[-0.200pt]{0.400pt}{4.818pt}}
\put(216.769,47.0){\makebox(0,0)[c]{4}}
\put(258.461,88.000){\rule[-0.200pt]{0.400pt}{4.818pt}}
\put(258.461,47.0){\makebox(0,0)[c]{5}}
\put(300.154,88.000){\rule[-0.200pt]{0.400pt}{4.818pt}}
\put(300.154,47.0){\makebox(0,0)[c]{6}}
\put(341.846,88.000){\rule[-0.200pt]{0.400pt}{4.818pt}}
\put(341.846,47.0){\makebox(0,0)[c]{7}}
\put(383.538,88.000){\rule[-0.200pt]{0.400pt}{4.818pt}}
\put(383.538,47.0){\makebox(0,0)[c]{8}}
\put(425.230,88.000){\rule[-0.200pt]{0.400pt}{4.818pt}}
\put(425.230,47.0){\makebox(0,0)[c]{9}}
\put(466.923,88.000){\rule[-0.200pt]{0.400pt}{4.818pt}}
\put(466.923,47.0){\makebox(0,0)[c]{10}}
\put(508.615,88.000){\rule[-0.200pt]{0.400pt}{4.818pt}}
\put(508.615,47.0){\makebox(0,0)[c]{11}}
\put(550.000,0.0){\makebox(0,0)[c]{$\log(B)$}}
\put(145.596,868.000){\rule[-0.200pt]{9.636pt}{0.400pt}}
\put(165.596,848.000){\rule[-0.200pt]{0.400pt}{9.636pt}}
\put(174.494,594.331){\rule[-0.200pt]{9.636pt}{0.400pt}}
\put(194.494,574.331){\rule[-0.200pt]{0.400pt}{9.636pt}}
\put(203.393,546.095){\rule[-0.200pt]{9.636pt}{0.400pt}}
\put(223.393,526.095){\rule[-0.200pt]{0.400pt}{9.636pt}}
\put(232.292,569.190){\rule[-0.200pt]{9.636pt}{0.400pt}}
\put(252.292,549.190){\rule[-0.200pt]{0.400pt}{9.636pt}}
\put(261.191,534.759){\rule[-0.200pt]{9.636pt}{0.400pt}}
\put(281.191,514.759){\rule[-0.200pt]{0.400pt}{9.636pt}}
\put(290.090,515.418){\rule[-0.200pt]{9.636pt}{0.400pt}}
\put(310.090,495.418){\rule[-0.200pt]{0.400pt}{9.636pt}}
\put(318.989,491.984){\rule[-0.200pt]{9.636pt}{0.400pt}}
\put(338.989,471.984){\rule[-0.200pt]{0.400pt}{9.636pt}}
\put(347.888,480.625){\rule[-0.200pt]{9.636pt}{0.400pt}}
\put(367.888,460.625){\rule[-0.200pt]{0.400pt}{9.636pt}}
\put(376.787,470.204){\rule[-0.200pt]{9.636pt}{0.400pt}}
\put(396.787,450.204){\rule[-0.200pt]{0.400pt}{9.636pt}}
\put(405.685,459.571){\rule[-0.200pt]{9.636pt}{0.400pt}}
\put(425.685,439.571){\rule[-0.200pt]{0.400pt}{9.636pt}}
\put(434.584,453.252){\rule[-0.200pt]{9.636pt}{0.400pt}}
\put(454.584,433.252){\rule[-0.200pt]{0.400pt}{9.636pt}}
\put(463.483,448.875){\rule[-0.200pt]{9.636pt}{0.400pt}}
\put(483.483,428.875){\rule[-0.200pt]{0.400pt}{9.636pt}}
\put(492.382,444.224){\rule[-0.200pt]{9.636pt}{0.400pt}}
\put(512.382,424.224){\rule[-0.200pt]{0.400pt}{9.636pt}}
\put(510.000,441.705){\rule[-0.200pt]{9.636pt}{0.400pt}}
\put(530.000,421.705){\rule[-0.200pt]{0.400pt}{9.636pt}}
\put(50.0,387.383){\rule[-0.200pt]{120.450pt}{0.400pt}}
\put(40.0,387.383){\makebox(0,0)[r]{$\thetaH(V)$}}
\end{picture}}
\hfil}
\hbox to\hsize{
\hfil$S_6$\hfil}}\]
\vfil
\penalty600
\[\vbox{%
\hsize\graphsize
\hbox to\hsize{\hfil
{\setlength{\unitlength}{0.240900pt}
\begin{picture}(600,900)(0,0)
\put(50.000,88.000){\rule[-0.200pt]{120.450pt}{0.400pt}}\put(50.000,88.000){\rule[-0.200pt]{0.400pt}{192.720pt}}\put(50.000,88.000){\rule[-0.200pt]{0.400pt}{4.818pt}}
\put(50.000,47.0){\makebox(0,0)[c]{0}}
\put(91.692,88.000){\rule[-0.200pt]{0.400pt}{4.818pt}}
\put(91.692,47.0){\makebox(0,0)[c]{1}}
\put(133.385,88.000){\rule[-0.200pt]{0.400pt}{4.818pt}}
\put(133.385,47.0){\makebox(0,0)[c]{2}}
\put(175.077,88.000){\rule[-0.200pt]{0.400pt}{4.818pt}}
\put(175.077,47.0){\makebox(0,0)[c]{3}}
\put(216.769,88.000){\rule[-0.200pt]{0.400pt}{4.818pt}}
\put(216.769,47.0){\makebox(0,0)[c]{4}}
\put(258.461,88.000){\rule[-0.200pt]{0.400pt}{4.818pt}}
\put(258.461,47.0){\makebox(0,0)[c]{5}}
\put(300.154,88.000){\rule[-0.200pt]{0.400pt}{4.818pt}}
\put(300.154,47.0){\makebox(0,0)[c]{6}}
\put(341.846,88.000){\rule[-0.200pt]{0.400pt}{4.818pt}}
\put(341.846,47.0){\makebox(0,0)[c]{7}}
\put(383.538,88.000){\rule[-0.200pt]{0.400pt}{4.818pt}}
\put(383.538,47.0){\makebox(0,0)[c]{8}}
\put(425.230,88.000){\rule[-0.200pt]{0.400pt}{4.818pt}}
\put(425.230,47.0){\makebox(0,0)[c]{9}}
\put(466.923,88.000){\rule[-0.200pt]{0.400pt}{4.818pt}}
\put(466.923,47.0){\makebox(0,0)[c]{10}}
\put(508.615,88.000){\rule[-0.200pt]{0.400pt}{4.818pt}}
\put(508.615,47.0){\makebox(0,0)[c]{11}}
\put(550.000,0.0){\makebox(0,0)[c]{$\log(B)$}}
\put(145.596,868.000){\rule[-0.200pt]{9.636pt}{0.400pt}}
\put(165.596,848.000){\rule[-0.200pt]{0.400pt}{9.636pt}}
\put(174.494,524.800){\rule[-0.200pt]{9.636pt}{0.400pt}}
\put(194.494,504.800){\rule[-0.200pt]{0.400pt}{9.636pt}}
\put(203.393,694.667){\rule[-0.200pt]{9.636pt}{0.400pt}}
\put(223.393,674.667){\rule[-0.200pt]{0.400pt}{9.636pt}}
\put(232.292,613.306){\rule[-0.200pt]{9.636pt}{0.400pt}}
\put(252.292,593.306){\rule[-0.200pt]{0.400pt}{9.636pt}}
\put(261.191,627.297){\rule[-0.200pt]{9.636pt}{0.400pt}}
\put(281.191,607.297){\rule[-0.200pt]{0.400pt}{9.636pt}}
\put(290.090,585.130){\rule[-0.200pt]{9.636pt}{0.400pt}}
\put(310.090,565.130){\rule[-0.200pt]{0.400pt}{9.636pt}}
\put(318.989,573.550){\rule[-0.200pt]{9.636pt}{0.400pt}}
\put(338.989,553.550){\rule[-0.200pt]{0.400pt}{9.636pt}}
\put(347.888,568.752){\rule[-0.200pt]{9.636pt}{0.400pt}}
\put(367.888,548.752){\rule[-0.200pt]{0.400pt}{9.636pt}}
\put(376.787,557.092){\rule[-0.200pt]{9.636pt}{0.400pt}}
\put(396.787,537.092){\rule[-0.200pt]{0.400pt}{9.636pt}}
\put(405.685,550.296){\rule[-0.200pt]{9.636pt}{0.400pt}}
\put(425.685,530.296){\rule[-0.200pt]{0.400pt}{9.636pt}}
\put(434.584,546.446){\rule[-0.200pt]{9.636pt}{0.400pt}}
\put(454.584,526.446){\rule[-0.200pt]{0.400pt}{9.636pt}}
\put(463.483,540.749){\rule[-0.200pt]{9.636pt}{0.400pt}}
\put(483.483,520.749){\rule[-0.200pt]{0.400pt}{9.636pt}}
\put(492.382,534.793){\rule[-0.200pt]{9.636pt}{0.400pt}}
\put(512.382,514.793){\rule[-0.200pt]{0.400pt}{9.636pt}}
\put(510.000,532.716){\rule[-0.200pt]{9.636pt}{0.400pt}}
\put(530.000,512.716){\rule[-0.200pt]{0.400pt}{9.636pt}}
\put(50.0,474.961){\rule[-0.200pt]{120.450pt}{0.400pt}}
\put(40.0,474.961){\makebox(0,0)[r]{$\thetaH(V)$}}
\end{picture}}
\hfil}
\hbox to\hsize{
\hfil$S_7$\hfil}}
\vbox{%
\hsize\graphsize
\hbox to\hsize{\hfil
{\setlength{\unitlength}{0.240900pt}
\begin{picture}(600,900)(0,0)
\put(50.000,88.000){\rule[-0.200pt]{120.450pt}{0.400pt}}\put(50.000,88.000){\rule[-0.200pt]{0.400pt}{192.720pt}}\put(50.000,88.000){\rule[-0.200pt]{0.400pt}{4.818pt}}
\put(50.000,47.0){\makebox(0,0)[c]{0}}
\put(91.692,88.000){\rule[-0.200pt]{0.400pt}{4.818pt}}
\put(91.692,47.0){\makebox(0,0)[c]{1}}
\put(133.385,88.000){\rule[-0.200pt]{0.400pt}{4.818pt}}
\put(133.385,47.0){\makebox(0,0)[c]{2}}
\put(175.077,88.000){\rule[-0.200pt]{0.400pt}{4.818pt}}
\put(175.077,47.0){\makebox(0,0)[c]{3}}
\put(216.769,88.000){\rule[-0.200pt]{0.400pt}{4.818pt}}
\put(216.769,47.0){\makebox(0,0)[c]{4}}
\put(258.461,88.000){\rule[-0.200pt]{0.400pt}{4.818pt}}
\put(258.461,47.0){\makebox(0,0)[c]{5}}
\put(300.154,88.000){\rule[-0.200pt]{0.400pt}{4.818pt}}
\put(300.154,47.0){\makebox(0,0)[c]{6}}
\put(341.846,88.000){\rule[-0.200pt]{0.400pt}{4.818pt}}
\put(341.846,47.0){\makebox(0,0)[c]{7}}
\put(383.538,88.000){\rule[-0.200pt]{0.400pt}{4.818pt}}
\put(383.538,47.0){\makebox(0,0)[c]{8}}
\put(425.230,88.000){\rule[-0.200pt]{0.400pt}{4.818pt}}
\put(425.230,47.0){\makebox(0,0)[c]{9}}
\put(466.923,88.000){\rule[-0.200pt]{0.400pt}{4.818pt}}
\put(466.923,47.0){\makebox(0,0)[c]{10}}
\put(508.615,88.000){\rule[-0.200pt]{0.400pt}{4.818pt}}
\put(508.615,47.0){\makebox(0,0)[c]{11}}
\put(550.000,0.0){\makebox(0,0)[c]{$\log(B)$}}
\put(145.596,868.000){\rule[-0.200pt]{9.636pt}{0.400pt}}
\put(165.596,848.000){\rule[-0.200pt]{0.400pt}{9.636pt}}
\put(174.494,587.200){\rule[-0.200pt]{9.636pt}{0.400pt}}
\put(194.494,567.200){\rule[-0.200pt]{0.400pt}{9.636pt}}
\put(203.393,478.000){\rule[-0.200pt]{9.636pt}{0.400pt}}
\put(223.393,458.000){\rule[-0.200pt]{0.400pt}{9.636pt}}
\put(232.292,479.137){\rule[-0.200pt]{9.636pt}{0.400pt}}
\put(252.292,459.137){\rule[-0.200pt]{0.400pt}{9.636pt}}
\put(261.191,447.531){\rule[-0.200pt]{9.636pt}{0.400pt}}
\put(281.191,427.531){\rule[-0.200pt]{0.400pt}{9.636pt}}
\put(290.090,403.103){\rule[-0.200pt]{9.636pt}{0.400pt}}
\put(310.090,383.103){\rule[-0.200pt]{0.400pt}{9.636pt}}
\put(318.989,383.620){\rule[-0.200pt]{9.636pt}{0.400pt}}
\put(338.989,363.620){\rule[-0.200pt]{0.400pt}{9.636pt}}
\put(347.888,365.556){\rule[-0.200pt]{9.636pt}{0.400pt}}
\put(367.888,345.556){\rule[-0.200pt]{0.400pt}{9.636pt}}
\put(376.787,350.849){\rule[-0.200pt]{9.636pt}{0.400pt}}
\put(396.787,330.849){\rule[-0.200pt]{0.400pt}{9.636pt}}
\put(405.685,338.251){\rule[-0.200pt]{9.636pt}{0.400pt}}
\put(425.685,318.251){\rule[-0.200pt]{0.400pt}{9.636pt}}
\put(414.000,334.727){\rule[-0.200pt]{9.636pt}{0.400pt}}
\put(434.000,314.727){\rule[-0.200pt]{0.400pt}{9.636pt}}
\put(434.584,328.410){\rule[-0.200pt]{9.636pt}{0.400pt}}
\put(454.584,308.410){\rule[-0.200pt]{0.400pt}{9.636pt}}
\put(459.804,319.996){\rule[-0.200pt]{9.636pt}{0.400pt}}
\put(479.804,299.996){\rule[-0.200pt]{0.400pt}{9.636pt}}
\put(510.000,308.383){\rule[-0.200pt]{9.636pt}{0.400pt}}
\put(530.000,288.383){\rule[-0.200pt]{0.400pt}{9.636pt}}
\put(50.0,223.880){\rule[-0.200pt]{120.450pt}{0.400pt}}
\put(40.0,223.880){\makebox(0,0)[r]{$\thetaH(V)$}}
\end{picture}}
\hfil}
\hbox to\hsize{
\hfil$S_8$\hfil}}\]
\vfil
\penalty 600
\vfilneg
We finish with tables of numerical results.
The value of $\thetaH^{\text{stat}}(V)$ is obtained
from the pairs $(B_i,\nUH(B_i))$ as described in section
\ref{section:stats}. We denote by $\zeta^*_{E_i}(1)$
the limit
\[\zeta^*_{E_i}(1)=\lim_{s\to 1}(s-1)^{t_i}\zeta_{E_i}^*(s),\]
where $t_i$ is the number of components of $E_i$.
Note that for the examples
with a Picard group of rank $2$, $C_2$ is equal to $1$.
\par
\ifx\FirstRKBcolumn\undefined
\createnewtable{RKB}{20}
\forRKB
\ifx\showallcolumns\undefined
\def\showallcolumns{}\fi
\gaddtomacro\showallcolumns{\showRKBcolumns}
\fi
\FirstRKBcolumn
{\text{Surface}}%
{B}%
{\nUH(B)}%
{\alpha(V)}%
{ad/cb}%
{\zeta^*_{E_1}(1)}%
{ab/cd}%
{\zeta^*_{E_2}(1)}%
{ac/bd}%
{\zeta^*_{E_3}(1)}%
{\lambda_{3}\omegaH(V(\mathbf Q_{3}))}%
{p_{0}}%
{\lambda_{p_{0}}\omegaH(V(\mathbf Q_{p_{0}}))}%
{C_0}%
{C_1}%
{C_3}%
{\omegaH(V(\mathbf R))}%
{\thetaH(V)}%
{\nUH(B)/\thetaH(V)B\log(B)}%
{\thetaH^{\mathrm{stat}}(V)/\thetaH(V)}%
\addRKBcolumn
{S_1}%
{20000}%
{75984}%
{2}%
{2}%
{8.146241\times10^{-1}}%
{1/8}%
{6.045998\times10^{-1}}%
{1/2}%
{8.146241\times10^{-1}}%
{4/9}%
{2}%
{3/8}%
{8.306815\times10^{-1}}%
{9.540383\times10^{-1}}%
{9.893865\times10^{-1}}%
{3.255161}%
{3.413500\times10^{-1}}%
{1.123839}%
{1.008178}%

\ifx\FirstRKBcolumn\undefined
\createnewtable{RKB}{20}
\forRKB
\ifx\showallcolumns\undefined
\def\showallcolumns{}\fi
\gaddtomacro\showallcolumns{\showRKBcolumns}
\fi
\FirstRKBcolumn
{\text{Surface}}%
{H}%
{\nUH(B)}%
{\alpha(V)}%
{ad/cb}%
{\zeta^*_{E_1}(1)}%
{ab/cd}%
{\zeta^*_{E_2}(1)}%
{ac/bd}%
{\zeta^*_{E_3}(1)}%
{\lambda'_{3}\omegaH(V(\mathbf Q_{3}))}%
{p_{0}}%
{\lambda'_{p_{0}}\omegaH(V(\mathbf Q_{p_{0}}))}%
{C_0}%
{C_1}%
{C_3}%
{\omegaH(V(\mathbf R))}%
{\thetaH(V)}%
{\nuH(B)/\thetaH(V)B\log(B)}%
{\thetaH^{\mathrm{stat}}(V)/\thetaH(V)}%
\addRKBcolumn
{S_2}%
{20000}%
{49608}%
{2}%
{5}%
{1.163730}%
{1/125}%
{6.045998\times10^{-1}}%
{1/5}%
{1.163730}%
{4/9}%
{5}%
{96/125}%
{3.493824\times10^{-1}}%
{8.704106\times10^{-1}}%
{9.906098\times10^{-1}}%
{1.360417}%
{2.290769\times10^{-1}}%
{1.093332}%
{0.958517}%

\ifx\FirstRKBcolumn\undefined
\createnewtable{RKB}{21}
\forRKB
\ifx\showallcolumns\undefined
\def\showallcolumns{}\fi
\gaddtomacro\showallcolumns{\showRKBcolumns}
\fi
\FirstRKBcolumn
{\text{Surface}}%
{H}%
{\nUH(B)}%
{\alpha(V)}%
{ad/cb}%
{\zeta^*_{E_1}(1)}%
{ab/cd}%
{\zeta^*_{E_2}(1)}%
{ac/bd}%
{\zeta^*_{E_3}(1)}%
{\lambda'_{3}\omegaH(V(\mathbf Q_{3}))}%
{p_{0}}%
{\lambda'_{p_{0}}\omegaH(V(\mathbf Q_{p_{0}}))}%
{C_0}%
{C_1}%
{C_3}%
{\omegaH(V(\mathbf R))}%
{\thetaH(V)}%
{\nUH(B)/\thetaH(V)B\log(B)}%
{\thetaH^{\mathrm{stat}}(V)/\thetaH(V)}%
\addRKBcolumn
{S_3}%
{20000}%
{78980}%
{2}%
{3}%
{1.017615}%
{1/27}%
{6.045998\times10^{-1}}%
{1/3}%
{1.017615}%
{4/9}%
{}%
{}%
{3.066383\times10^{-1}}%
{9.762028\times10^{-1}}%
{9.892790\times10^{-1}}%
{2.221359}%
{3.660885\times10^{-1}}%
{1.089213}%
{1.021637}%

\showRKBcolumns
\par
\vfil
\penalty -600
\noindent
For the examples
with a Picard group of rank $3$, $C_3$ is equal to $1$.
\par
\ifx\FirstRKCcolumn\undefined
\createnewtable{RKC}{20}
\forRKC
\ifx\showallcolumns\undefined
\def\showallcolumns{}\fi
\gaddtomacro\showallcolumns{\showRKCcolumns}
\fi
\FirstRKCcolumn
{\text{Surface}}%
{B}%
{\nUH(B)}%
{\alpha(V)}%
{ad/bc}%
{\zeta^*_{E_1}(1)}%
{ab/cd}%
{\zeta^*_{E_2}(1)}%
{ac/bd}%
{\zeta^*_{E_3}(1)}%
{\lambda_{3}\omegaH(V(\mathbf Q_{3}))}%
{p_{0}}%
{\lambda_{p_{0}}\omegaH(V(\mathbf Q_{p_{0}}))}%
{C_0}%
{C_1}%
{C_2}%
{\omegaH(V(\mathbf R))}%
{\thetaH(V)}%
{\nUH(B)/\thetaH(V)B\log(B)^{2}}%
{\thetaH^{\mathrm{stat}}(V)/\thetaH(V)}%
\addRKCcolumn
{S_4}%
{100000}%
{3051198}%
{1}%
{1}%
{6.045998\times10^{-1}}%
{1/4}%
{8.146241\times10^{-1}}%
{1}%
{6.045998\times10^{-1}}%
{16/27}%
{2}%
{27/64}%
{8.306815\times10^{-1}}%
{9.540383\times10^{-1}}%
{7.827314\times10^{-1}}%
{4.105301}%
{1.895795\times10^{-1}}%
{1.214249}%
{0.980969}%

\ifx\FirstRKCcolumn\undefined
\createnewtable{RKC}{20}
\forRKC
\ifx\showallcolumns\undefined
\def\showallcolumns{}\fi
\gaddtomacro\showallcolumns{\showRKCcolumns}
\fi
\FirstRKCcolumn
{\text{Surface}}%
{H}%
{\nUH(B)}%
{\alpha(V)}%
{ad/bc}%
{\zeta^*_{E_1}(1)}%
{ab/cd}%
{\zeta^*_{E_2}(1)}%
{ac/bd}%
{\zeta^*_{E_3}(1)}%
{\lambda'_{3}\omegaH(V(\mathbf Q_{3}))}%
{p_{0}}%
{\lambda'_{p_{0}}\omegaH(V(\mathbf Q_{p_{0}}))}%
{C_0}%
{C_1}%
{C_2}%
{\omegaH(V(\mathbf R))}%
{\thetaH(V)}%
{\nUH(B)/\thetaH(V)B\log(B)^{2}}%
{\thetaH^{\mathrm{stat}}(V)/\thetaH(V)}%
\addRKCcolumn
{S_5}%
{100000}%
{1976482}%
{1}%
{1}%
{6.045998\times10^{-1}}%
{1/25}%
{1.163730}%
{1}%
{6.045998\times10^{-1}}%
{16/27}%
{5}%
{13824/15625}%
{3.493824\times10^{-1}}%
{8.704106\times10^{-1}}%
{8.112747\times10^{-1}}%
{2.347970}%
{1.291945\times10^{-1}}%
{1.154191}%
{1.056178}%

\ifx\FirstRKCcolumn\undefined
\createnewtable{RKC}{20}
\forRKC
\ifx\showallcolumns\undefined
\def\showallcolumns{}\fi
\gaddtomacro\showallcolumns{\showRKCcolumns}
\fi
\FirstRKCcolumn
{\text{Surface}}%
{H}%
{\nUH(B)}%
{\alpha(V)}%
{ad/bc}%
{\zeta^*_{E_1}(1)}%
{ab/cd}%
{\zeta^*_{E_2}(1)}%
{ac/bd}%
{\zeta^*_{E_3}(1)}%
{\lambda'_{3}\omegaH(V(\mathbf Q_{3}))}%
{p_{0}}%
{\lambda'_{p_{0}}\omegaH(V(\mathbf Q_{p_{0}}))}%
{C_0}%
{C_1}%
{C_2}%
{\omegaH(V(\mathbf R))}%
{\thetaH(V)}%
{n_\greeksubscript{U,\mathbf H}(B)/\thetaH(V)B\log(B)^{2}}%
{\thetaH^{\mathrm{stat}}(V)/\thetaH(V)}%
\addRKCcolumn
{S_6}%
{100000}%
{3420784}%
{1}%
{1}%
{6.045998\times10^{-1}}%
{1/49}%
{1.265025}%
{1}%
{6.045998\times10^{-1}}%
{16/27}%
{7}%
{186624/117649}%
{3.066383\times10^{-1}}%
{9.297617\times10^{-1}}%
{9.228033\times10^{-1}}%
{1.910125}%
{2.184437\times10^{-1}}%
{1.181448}%
{1.044971}%

\ifx\FirstRKCcolumn\undefined
\createnewtable{RKC}{20}
\forRKC
\ifx\showallcolumns\undefined
\def\showallcolumns{}\fi
\gaddtomacro\showallcolumns{\showRKCcolumns}
\fi
\FirstRKCcolumn
{\text{Surface}}%
{H}%
{\nUH(B)}%
{\alpha(V)}%
{ad/bc}%
{\zeta^*_{E_1}(1)}%
{ab/cd}%
{\zeta^*_{E_2}(1)}%
{ac/bd}%
{\zeta^*_{E_3}(1)}%
{\lambda'_{3}\omegaH(V(\mathbf Q_{3}))}%
{p_{0}}%
{\lambda'_{p_{0}}\omegaH(V(\mathbf Q_{p_{0}}))}%
{C_0}%
{C_1}%
{C_2}%
{\omegaH(V(\mathbf R))}%
{\thetaH(V)}%
{n_\greeksubscript{U,\mathbf H}(B)/\thetaH(V)B\log(B)^{2}}%
{\thetaH^{\mathrm{stat}}(V)/\thetaH(V)}%
\addRKCcolumn
{S_7}%
{100000}%
{1966160}%
{1}%
{1}%
{6.045998\times10^{-1}}%
{4/9}%
{1.028996}%
{1}%
{6.045998\times10^{-1}}%
{16/27}%
{2}%
{27/64}%
{8.306815\times10^{-1}}%
{8.196347\times10^{-1}}%
{8.294515\times10^{-1}}%
{2.430506}%
{1.290720\times10^{-1}}%
{1.149252}%
{0.981831}%

\showRKCcolumns
\par
\vfil
\penalty -600
\noindent
For the last example we have $C_2=C_3=1$ and $E_1=E_2=E_3$
and we get
\par
\ifx\FirstRKDcolumn\undefined
\createnewtable{RKD}{12}
\forRKD
\ifx\showallcolumns\undefined
\def\showallcolumns{}\fi
\gaddtomacro\showallcolumns{\showRKDcolumns}
\fi
\FirstRKDcolumn
{\text{Surface}}%
{B}%
{\nUH(B)}%
{\alpha(V)}%
{\zeta^*_{E_i}(1)}%
{\lambda_{3}\omegaH(V(\mathbf Q_{3}))}%
{C_0}%
{C_1}%
{\omegaH(V(\mathbf R))}%
{\thetaH(V)}%
{\nUH(B)/\thetaH(V)B\log(B)^{3}}%
{\thetaH^{\mathrm{stat}}(V)/\thetaH(V)}%
\addRKDcolumn
{S_8}%
{100000}%
{12137664}%
{7/18}%
{6.045998\times10^{-1}}%
{16/27}%
{3.066383\times10^{-1}}%
{5.129319\times10^{-1}}%
{6.121864}%
{4.904057\times10^{-2}}%
{1.621894}%
{1.024630}%

\showRKDcolumns
\par

\ifx\undefined\bysame
\newcommand{\bysame}{\leavevmode\hbox to3em{\hrulefill}\,}
\fi
\ifx\undefined\numero
\newcommand{\numero}{$\hbox{n}^{\hbox{\scriptsize o}}\,$}
\else\renewcommand{\numero}{$\hbox{n}^{\hbox{\scriptsize o}}\,$}
\fi


\end{document}